
\input amssym.def
\input amssym
\input xypic

\magnification = \magstep1
\hsize=14.4truecm
\vsize= 24.5truecm
\hoffset=0.4 truecm
\baselineskip=12pt

\font\big=cmbx10 scaled \magstep2
\font\sbig=cmbx10

\font\fascio=eusm10

\font\titol=cmcsc8
\font\proj=msbm10 scaled \magstep2
\font\eightpt=cmr8

\def\smallskip{\vskip 0.2truecm}
\def\medskip{\vskip 0.5truecm}
\def\bigskip{\vskip 0.7truecm}

\def\tit{\medskip \noindent} 

\def\proof{\par\smallskip\noindent{\it Proof.}\ }
\def\endproof{\hfill$\square$\smallskip}
\definemorphism{mapsin}\solid\tip\ahook

\def\sp#1{\Bbb P^{#1}}
\def\cO{{\cal O}}

\def\ga{p_{\scriptscriptstyle a}}
\def\cc{c_{\scriptscriptstyle 1}}

\headline={\ifnum\pageno=1\hfil\else
{\ifodd\pageno
\rightheadline\else\leftheadline\fi}\fi}
\def\rightheadline{\titol\hfil 
On smooth rational threefolds of $\sp 5 \ \ldots$
\hfil{\tenrm\folio}}                    
\def\leftheadline{{\tenrm\folio}\titol\hfil 
Mezzetti - Portelli \hfil}

\nopagenumbers

{\big 
\hbox{On smooth rational threefolds of 
{\proj P}${}^{\hbox{\sbig 5}}$ with}}

\medskip

{\big 
\hbox{rational non-special hyperplane section}}
\footnote{}{1991 {\it Mathematics Subject Classification.}
14C20, 14J30, 14M07, 14M20}
\footnote{}{{\it Keywords and phrases.} Codimension two, rational
threefolds, birational maps.}
\bigskip

\noindent{\it Emilia Mezzetti \ } and {\ \it 
Dario Portelli}
\medskip

{\eightpt {
\noindent Dipartimento di 
Scienze Matematiche, \ Universit\`a di Trieste,

\noindent Piazzale Europa 1,
 
\noindent 34127 Trieste, ITALY

\noindent e-mail: mezzette@univ.trieste.it, 
porteda@univ.trieste.it}}

\bigskip

\tit
{\bf Introduction}
\tit

Linear systems of plane curves defining rational 
surfaces of $\sp 4$
have been recently extensively studied. In particular:
Alexander ([A], [A1]) gave  theorems
of existence and unicity for linear systems defining 
smooth surfaces
in $\sp 4$ with non-special or speciality one 
hyperplane sections,
Catanese-Franciosi ([CF]) and Catanese-Hulek ([CH]) 
have performed a fine
classification of base loci of linear systems 
defining smooth surfaces.

But nothing has been written in recent years
about linear systems of surfaces of $\sp 3$ 
giving rational threefolds
in $\sp r$. Nevertheless, something on the 
subject may be found in some classical papers by 
F.\thinspace Jongmans and U.\thinspace Morin 
([J],[M],[M1]).

Here we study, from the above point of view, 
the smooth threefolds of 
$\sp 5$ having a rational non-special surface 
of $\sp 4$ as general 
hyperplane section. It is known ([I], [I1]) that 
there are exactly five of such threefolds: 
$\sp 2\times\sp 1,$ Del Pezzo, Castelnuovo and
the scrolls of Bordiga and Palatini (of degrees 
3,4,...,7 respectively). Among them we find 
all smooth threefolds which are scrolls over a 
rational surface ([O]). It has  been in fact the 
wish of better understanding the geometry of the 
scrolls of Bordiga and Palatini which originally 
motivated our interest to this subject.

Our main result is that, for degree 3,4,5,6 the
threefold $X$ contains a line $L$ such that the
projection of center $L$ gives a birational map
between $X$ and $\sp 3.$
In this situation
the linear system $\vert \Sigma\vert$ of surfaces in $\sp 3$
defining the inverse map of the projection
is particularly simple; in fact,
if $deg \  X=d$, then $deg \ \Sigma =d-1$; moreover 
the restriction of $\vert \Sigma\vert$ to a general plane is
a linear system giving a general hyperplane section of $X$ 
containing $L.$

We completely describe the base locus $B$ of
$\vert\Sigma\vert$ in the four cases; precisely we
determine the degree and the arithmetic genus
for the irreducible components of $B$, the infinitesimal
structure we need to consider on them, and the way 
these components reciprocally intersect. 
An interesting consequence of this construction is
the description of $X$ as 
a suitable blowing-down of the blowing-up of $\sp 3$
along $B.$

The situation is different for the Palatini scroll, 
because a similar line never
exists. In a forthcoming paper ([MP]) we will give
a different construction of a birational map
$X\dashrightarrow\sp 3$ in this last case.

\smallskip
\noindent
The results on $\vert\Sigma\vert$ we obtain are 
collected in the following table (where \lq\lq genus"
always means \lq\lq arithmetic genus").
 \smallskip

\halign{\hfil\bf# & \quad\hfil#\hfil & \quad\hfil#\hfil & \quad#\hfil    \cr
 Variety & {\bf Degree} & $deg\thinspace\Sigma$ & {\bf Base locus}  \cr
 
&&& \cr
$\sp 2\times\sp 1$ & 3 & 2 &  $\bullet$\ a point $P$                         \cr
                   &   &   &  $\bullet$\ a line                                        \cr
                   &   &   &                                                           \cr
Del Pezzo          & 4 & 3 &  $\bullet$\ a quintic curve of genus $2$                  \cr 
                   &   &   &                                                           \cr
Castelnuovo        & 5 & 4 &  $\bullet$\ a  curve $B_1$ of degree $7$ and   genus $3$                 \cr
                   &   &   &  \ \ \thinspace  with a $5$-secant line $B_2$    \cr
                   &   &   &  $\bullet$\ the first infinitesimal nbh. of $B_2$         \cr
                   &   &   &                                                           \cr 
Bordiga scroll     & 6 & 5 &  $\bullet$\ the first infinitesimal nbh. of               \cr
																			&   &   &  \ \ \thinspace a curve $B_1$ of degree $3$ and genus $0$   \cr
                   &   &   &  $\bullet$\ a curve $B_2$ of degree  $7$ and  genus $0$      \cr
 																		&   &   &  \ \ \thinspace with $deg\ (B_1\cap B_2)=12$ \cr}

It is natural to ask whether we can solve the inverse
problem: characterize, among the linear systems satisfying
the conditions of the above table, those defining 
Del Pezzo, Castelnuovo and Bordiga threefolds. 
We are able to completely solve this problem: the answer
is that the base curve $B$ has to be locally 
Cohen--Macaulay and, roughly speaking, its general plane section 
has to satisfy the conditions given by Catanese--Franciosi.

Let us say a few words about the methods employed 
and the organization of the paper. 

The study  of the threefolds $X$ of degree $3$ 
and $4,$ which is briefly reviewed in \S\thinspace 1, 
is classical and easy. It motivates the analysis of
the general properties of birational maps which are 
projections from an embedded line $L\subset X,$ 
which is the content of \S\thinspace 2.

In the cases of  Castelnuovo $3$-folds 
and Bordiga scrolls
it is not  a priori obvious that a similar
projection exists. In the first case, a
rather deep  investigation
of the geometry of $X$, inspired from
[M], shows that in fact a suitable line always 
exists. We do this in \S\thinspace 3. 
 
We collect in \S\thinspace 4 the geometric properties
of the scrolls of Bordiga and Palatini we will need
in the subsequent sections.

If $X$ is a Bordiga scroll, the existence of a suitable 
line is proved after a detailed study of the family of
the adjunction maps
for the hyperplane sections of $X$ 
(\S\thinspace 5). Each of them
contracts $10$ lines to $10$ points of $\sp 2$; 
using a connectedness theorem of Debarre, 
we prove that for some hyperplane section $S$
these points are in special position, i.e. $7$ of
them are on a conic. The inverse image of such
conic on $S$ is the desired line on $X.$

In \S\thinspace 6 we show that for every Palatini
scroll a line as above does not exist.

In \S\thinspace 7 we reverse the point of view:
we start from a linear system $\vert\Sigma\vert$ of
surfaces in $\sp 3$ and determine necessary and sufficient 
conditions on the base locus in order that $\vert\Sigma
\vert$ defines a smooth threefold in $\sp 3$ 
of degree $3,4,5$ and $6$ respectively. One of the
key tool here is the numerical character of a curve. 

\smallskip

We will always work over an algebraically closed 
field $k$ of characteristic zero. 
We will assume throughout this paper that $X\subset
\sp 5$ is an irreducible, smooth, non degenerate 
threefold.
If $E\subset\sp n,$ we shall denote by $\langle 
E\rangle$ the linear span of $E.$

\medskip
{\bf  Acknowledgements.}

{\it This work has been done in the framework of the 
activities of Europroj. Both authors have been
supported by funds of MURST and  
GNSAGA, Progetto Strategico \lq\lq\thinspace
Applicazioni della Matematica per la Tecnologia
e la Societ\` a", Sottoprogetto \lq\lq\thinspace
Calcolo simbolico".}


\tit
{\bf  1.-\ Preliminaries.}
\smallskip


In this section we collect some classically known facts
about birational maps and describe from our point of view 
the two first examples of rational threefolds, namely
$\sp 2\times\sp 1$ and the Del Pezzo threefold, 
complete intersection of two quadric hypersurfaces.
In both cases a birational map to $\sp 3$ is simply given by
the projection from a suitable line contained in $X$.
At the end of the section,
for the reader convenience, we will state
a result of Catanese-Franciosi
([CF]) giving necessary and sufficient conditions
on a zero-dimensional subscheme $W\subset\sp 2$ 
for the very ampleness on the blow-up of  
$\sp 2$ along $W$ of the pull-back of a suitable 
linear system $\vert H\vert$ of plane curves. 

We recall first some classical facts on
birational maps.

\medskip
\noindent {\bf 1.1} \ {\it Fundamental points and exceptional 
divisors for a birational map.}

\smallskip

Let $f:X\dashrightarrow Y$ be a birational map between 
smooth varieties. A point $x\in X$
is called {\sl fundamental} for $f$ in the case $f$ is 
not regular at $x.$
Since $X$ is smooth (\lq\lq\thinspace normal" is sufficient), 
the set $F$ of the fundamental points for $f$ is a closed 
subset of $X,$ of codimension at least $2.$ We will call $F$
the {\sl fundamental locus} for $f.$ An irreducible 
subvariety $Z\subset X$ contained in the fundamental
locus will be also called {\sl fundamental}.

\smallskip 
\noindent
{\bf Van der Waerden'\thinspace s Purity Theorem.}\ {\sl 
Let $f:X\dashrightarrow Y$ be a birational map of smooth 
varieties, and let
$W\subset Y$ be a fundamental variety for $f^{-1}.$
Assume that $f(X)\cap W$ is dense in $W.$ Then any component $E$ of $f^{-1}(W)$ is of 
codimension $1$ in $X.$} 

\noindent For a proof we can refer to [EGA] (21.1), 
where the theorem is given under weaker assumptions.

\smallskip
 
We will call any $E\subset X$ as above an {\sl exceptional 
divisor} for $f.$ 

\smallskip 
\noindent
{\bf Lemma 1.1}\ {\sl  Let $g:\sp n\dashrightarrow X$ be a
birational map defined by the linear 
system $\vert\Sigma
\vert$ of hypersurfaces in $\sp n.$ 
Then any irreducible
component of the fundamental locus for $g$ 
is also a
component of $Bs(\vert\Sigma\vert ),$ 
and conversely.}
\proof
See [Z], Thm.\thinspace 15.

\medskip
Let us finally recall that a {\sl characteristic
curve} $\Gamma$ of a linear system 
$\vert\Sigma\vert$ of surfaces of $\sp 3$
is the free intersection of two general surfaces
of $\vert\Sigma\vert$.
If $E$ is an exceptional surface for $g^{-1}$ (notations are as
in the lemma above), and $B$ is the corresponding
base curve, then $deg \ E=deg \ (B\cap\Gamma)$.

\medskip
\noindent {\bf 1.2.}{\it The rational normal scroll}  $\sp 2\times\sp 1$.

\smallskip

Let $X\subset\sp 5$ be the image of $\sp 2\times\sp 1$
embedded in $\sp 5$ by the Segre map $s$. 
For $x$ a point of $\sp 1,$ we will denote 
 $F_1 = s(\sp 2\times x),$ a plane on $X.$
Moreover, let $l$ be a line in $\sp 2$ and
set $F_2=s(l\times\sp 1).$ $F_2$ is a quadric surface on $X.$ 
Finally, set $L:=F_1\cap F_2.$

The properties collected in the next 
proposition are classical ([SR]).

\smallskip 
\noindent
{\bf Proposition 1.2.1} {\sl 
\item{(i)} The projection map $ 
\pi_L :X \dashrightarrow  \sp 3 $ 
with center $L$ is birational;

\item{(ii)} the exceptional divisors of $\pi_L$ are
$F_1$ and $F_2,$ in particular $\pi_L(F_1)$ is a 
 single point $P$ of $\sp 3,$ and $\pi_L(F_2)$ 
is a line $B\subset\sp 3$ such that $P\notin B;$
 
\item{(iii)} the map $\pi_L^{-1}$ is defined by the 
linear system  of quadrics $\vert\Sigma\vert$ with 
 base  locus $B\cup P;$

\item{(iv)} the only exceptional divisor for $\pi_L^{-1}$ 
is the plane $\Phi =\langle B\cup P\rangle ,$ which is 
 contracted to the line $L .$
 
}

\proof A general hyperplane through $L$ cuts $X$ 
into a rational normal
scroll $S.$ Another general hyperplane through $L$ 
cuts $S$ into the
line $L$ plus a conic having only one point in common 
with $L.$ Finally,
a third general hyperplane through $L$ cuts this 
curve outside $L$ in
exactly one point. Therefore, for a general point 
$A\in X,$ the plane
$\langle L\cup A\rangle$ cuts $X$ outside $L$ only 
in $A,$ and (i) is proved.

The first part of (ii) is trivial and the second  
follows easily by 
remarking that $\pi_L$ is constant on any line 
of the ruling
on $F_2$ which does not contain $L.$

A general hyperplane section $H$ of $X$
intersects $L$ at one point, hence it projects
on a quadric $\Sigma$. $H$ meets $F_1$
and all lines of $F_2$, therefore
the base locus of
$\vert\Sigma\vert$ contains $B\cup P.$
Conversely, 
the linear system of the quadric surfaces in $\sp 3$ 
through $B\cup P$
has dimension $5,$ and (iii) follows.

(iv): let $R$ be a line of $\Phi$ through $P$; 
$R$ intersects any surface of $\vert \Sigma\vert$ in two
base points, so it contracts to a point of $L$.
\endproof

\noindent
{\bf Remark 1.2.2.} If we project from a line
$s(y\times\sp 1)$, $y\in\sp 2$, 
of the ruling on $X,$
then the projection is not birational, 
as it is easily seen.

\noindent
{\bf Remark 1.2.3.} The linear system $\vert\Sigma\vert$ 
contains a subsystem of dimension $3$ formed
by reducible quadrics, which are the union of the plane
$\Phi$ and of a variable plane. The corresponding 
hyperplane 
sections of $X$ are precisely those which are cut out
by hyperplanes containing $L$.

\noindent
{\bf Construction 1.2.4.}\ 
To better understand the projection $\pi_L$, 
let us perform the following construction.
Let $\sigma: \tilde\sp 3:= \tilde\sp 3(B\cup P)
\longrightarrow\sp 3$
be the blowing-up of $\sp 3$ along $B\cup P$. 
The linear system 
$\vert\Sigma\vert$ induces on $\tilde\sp 3$ 
the linear system $\Psi=\vert2H-E_B-E_P\vert$, where $H$
is the pull-back via $\sigma$ of the hyperplane divisor of
$\sp 3$, and $E_B, E_P$ are the exceptional divisors
over $B$ and $P$ respectively. The map defined by $\Psi$
is regular and fits in a commutative diagram:
$$
\diagram
\tilde\sp 3 \dto_\sigma \drto{} & \cr
\sp 3  \rdashed_{\pi_L^{-1}}|>\tip  & X 
\enddiagram
$$
But $\Psi$ is not very ample, in fact on
 $X$ the images of $E_B$ and $E_P$ are
the plane $F_1$ and the quadric $F_2$ which 
intersect along $L$. 

Let us consider 
another linear system on
$\tilde\sp 3$: $\Psi '= \vert 3H-E_B-E_P\vert$. It results 
to be very ample because the homogeneous 
ideal of $B\cup P$ in $\sp 3$ is generated in degree $2$
(see [Co],[BS]). Therefore $\Psi '$ defines an embedding 
$\tilde\sp 3 \longrightarrow \sp {14}$ with image $Y$, 
a smooth threefold of degree $23$.

Let now $\widetilde X$ be the blowing-up of $X$ along $L$.
The homogeneous ideal of $L$ in $X$ is clearly generated
in degree one; so by the  results quoted above,
the linear system  $\vert 2H_X-L\vert$ on $\widetilde X$
 is very ample and defines an embedding 
$\widetilde X\longrightarrow \sp {14}$. It is immediate
 to verify that the image is the variety $Y$. We have 
finally the following commutative diagram clarifying 
the geometry of the projection $\pi _L$, which results 
to be the composition of a blowing-up with
the inverse of a blowing-up:
$$
\diagram
 Y\rto|<<\ahook  &\sp {14}  \\
\tilde\sp 3\dto_{\sigma}  \rto^{\sim} \uto & \tilde X\dto 
\ulto|>\tip\\
\sp 3  \rdashed_{\pi_L^{-1}}|>\tip  & X
\enddiagram
$$

\noindent
{\bf Remark 1.2.5.}\ The planes $\sp 2\times a$ on $\sp 2\times\sp 1
\subset\sp 5$ correspond via $\pi_L$ to the planes in $\sp 3$ containing 
the line $B.$ The  quadrics
$R\times\sp 1\subset\sp 2\times\sp 1,$ where $R\subset\sp 2$ is a line,
correspond via $\pi_L$ to the planes in $\sp 3$ containing 
the point $P.$

\medskip
\noindent {\bf 1.3.} \ {\it Del Pezzo Threefold.}

\smallskip

Let $X\subset\sp 5$ be a Del Pezzo threefold, 
namely the complete intersection
of two quadric hypersurfaces in $\sp 5.$ Let 
$L\subset X$ be
any line. 
The following properties, similar to those of
Prop.1.2.1, hold  ([SR]):

\smallskip 
\noindent
{\bf Proposition 1.3.1} {\sl 
\item{(i)} The projection map $\pi_L :X \dashrightarrow  
\sp 3 $ with center $L$ is birational;

\item{(ii)} the map $\pi_L^{-1}$ is defined by a 
linear system  of cubic surfaces $\vert\Sigma\vert$ with 
base 
locus a general quintic curve $B$ 
of genus two;

\item{(iii)} the only exceptional divisor for $\pi_L^{-1}$ 
is the quadric $\Phi$ containing $B,$ which 
is contracted by $\pi^{-1}_L$  to $L;$

\item{(iv)} the lines contained in $X$ and meeting 
$L$ generate a singular ruled surface 
$T,$ of  degree $8,$ having $L$ as 
triple locus, which is the only exceptional divisor 
of $\pi_L .$ The surface $T$ is
contracted by $\pi_L$ to $B.$

}

\proof (i) is similar to (i) of Proposition 1.2.1: now 
a general section of $X$ with a $\sp 3$ containing $L$
is a quartic curve which splits in the union of $L$ with a
rational cubic, having $L$ as a chord. Hence a general plane 
through $L$ intersects $X$ at a unique point outside $L$.

A general hyperplane section $H$ of $X$
intersects $L$ at one point, hence it projects
on a cubic $\Sigma$. A general curve section $C$ of $X$
doesn't intersect $L$, so the characteristic curves of
$\vert\Sigma\vert$ are elliptic quartics. The linked
curve to such a quartic $D$, in the complete intersection
of two cubic surfaces containing it, is a quintic curve
$B$ of genus $2,$
meeting $D$ in $8$ points. Conversely, an easy
calculation with cohomology shows that
the linear system of the cubic surfaces in $\sp 3$ 
through $B$ has dimension $5$ and defines a rational map
whose image is a threefold of degree $4$.

$B$ is contained in a unique quadric $\Phi$, and is a 
divisor of type $(2,3)$ on $\Phi$. Clearly the lines of 
$\Phi$ which meet $B$ at three points are contracted by 
$\vert\Sigma\vert$, those meeting $B$ at two points
go to $L$. 

The points of the base curve $B$ of $\vert\Sigma\vert$
come from lines in $X$ meeting $L$; the degree of 
the surface $T$ they generate is clearly 
$deg \ (B\cap D)=8$.
Finally, $L$ is triple for $T$ because each point of $L$
comes from a trisecant line for $B.$
\endproof

\smallskip
\noindent
{\bf Remark 1.3.2.} As in the case of $\sp 2
\times \sp 1$, the linear system $\vert\Sigma\vert$ 
contains a subsystem of dimension $3$ formed
by reducible cubic surfaces, which are the union of the 
quadric $\Phi$ and of a variable plane. This subsystem
corresponds to hyperplane sections 
of $X$ containing $L$.

\smallskip
\noindent
{\bf Remark 1.3.3.} As in 1.2.4, we may factorize the rational 
map $\pi_L^{-1}$ through a variety $Y$ which is 
 isomorphic to both
the blowing-up of $\sp 3$ along the quintic $B$ 
and to the blowing-up
of $X$ along the line $L$. $Y$ is a 
threefold of degree $26$ in 
$\sp {15}$ and can be seen as the image of 
 the rational map defined by the linear system
of quartic surfaces 
$\vert 4H_{\scriptscriptstyle {\sp 3}}-B\vert$ of $\sp 3;$ 
alternatively, it is 
defined by the linear system $\vert 2H_X-L\vert$ on $X.$
This follows as in 1.2.4 from the known fact that
the homogeneous ideal of $B$ is generated in degree $3.$

\smallskip
\noindent
{\bf Remark 1.3.4.} The lines on $X$ such
that the projection with center one of them
gives a birational map from $X$ to $\sp 3$
(\lq\lq good centers of projection") form a 
family of dimension $3$ and $2,$ respectively
for $\sp 2\times\sp 1$ and for a Del Pezzo threefold.

\medskip
\noindent {\bf 1.4.} \ {\it The Catanese-Franciosi
conditions.}

\smallskip
Let $W$ denote a zero-dimensional subscheme of $\sp 2$
and let $\vert H\vert$ be a linear
system of plane curves. We recall here under what 
conditions on $W$ the pull-back
of $H$ on the blowing-up $\pi :S\to \sp 2$ of $\sp 2$
 along $W$ is very ample.
\smallskip

\noindent
(CF1) \ $deg\ W=5$ \ ($S$ is a Del Pezzo surface) 

\vskip -0.2 cm

\noindent
------------------------------------------------------------

\noindent
$W=\{ x_1,\ldots ,x_5\}$ and
$\vert H\vert =\vert 3H_{\scriptscriptstyle{\sp 2}}
-W\vert$ (here and in all other cases
$H_{\scriptscriptstyle{\sp 2}}$
denotes a line in $\sp 2$); then $\pi ^{*}H$ is
very ample if and only if:

\smallskip
\item{(0)} $W$ has no 
infinitely near points ;

\item{(1)} $h^0(H_{\scriptscriptstyle{\sp 2}}
-x_i-x_j-x_k)=0$, for any 
$i,j,k\in\lbrace 1,\dots ,5\rbrace$ .

\medskip

\noindent
(CF2) \ $deg\ W=8$ \ ($X$ is a Castelnuovo surface) 

\vskip -0.2 cm

\noindent
-----------------------------------------------------------------

\smallskip

\noindent $W=
\{ x_1,y_2\ldots ,y_8\}$ and $\vert H\vert =
\vert 4H_{\scriptscriptstyle{\sp 2}}
-2x_1-\sum_{\scriptscriptstyle {
2\leqslant i\leqslant 8}}y_i\vert$,
then $\pi ^{*}H$ is
very ample if and only if:

\smallskip
\item{(0)} at most one among $y_2,\dots,y_8$ is
infinitely near to $x_1$ ;

\item{(1)} $h^0(H_{\scriptscriptstyle{\sp 2}}
-x_1-y_i-y_j)=0$, for any 
$i,j\in\lbrace 2,\dots ,8\rbrace$ ;

\item{(2)} $h^0(H_{\scriptscriptstyle{\sp 2}}
-\sum_{i\in\Delta} y_i)=0$ for 
$\#\Delta\geqslant 4$ ;

\item{(3)} $h^0(2H_{\scriptscriptstyle{\sp 2}}
-x_1-\sum_{i\neq j}y_j)=0$, for  
any $j\in\lbrace 2,\dots ,8\rbrace$ .

\medskip

\noindent
(CF3) \ $deg\ W=10$ \ ($S$ is a Bordiga surface) 

\vskip -0.2 cm

\noindent
------------------------------------------------------------ 

\noindent $W=\{ y_1,\ldots ,y_{10}\}$, and
$\vert H\vert =\vert 4H_{\scriptscriptstyle{\sp 2}}
-W\vert$; then $\pi ^{*}H$ is
very ample if and only if:
\smallskip 
\item{(0)} $W$ has no 
infinitely near points ;

\item{(1)} $h^0(H_{\scriptscriptstyle{\sp 2}}
-\sum_{i\in\Delta} y_i)=0$, for 
$\#\Delta\geqslant 4$ ;

\item{(2)} $h^0(2H_{\scriptscriptstyle{\sp 2}}
-\sum_{i\in\Delta} y_i)=0$, for 
$\#\Delta\geqslant 8$ ;

\item{(3)} $h^0(3H_{\scriptscriptstyle{\sp 2}}
-\sum_{\scriptscriptstyle {
1\leqslant i\leqslant 10}} y_i)=0.$

\smallskip\noindent
Actually, we will need to consider also another
linear system $\vert H'\vert$ of plane curves
defining the Bordiga surface; precisely, $\vert H'\vert$
is obtained from $\vert H\vert$ by means of the
standard quadratic transformation centered at three
non collinear points among the $y_i$'s. Then the
above conditions modify to (notations are fresh):

\smallskip

\noindent
(CF3') \ $deg\ W=10$ \ ($S$ is again a Bordiga surface) 

\vskip -0.2 cm

\noindent
--------------------------------------------------------------------

\noindent $W=\{ x_1,x_2,x_3,y_4,\ldots 
,y_{10}\}$, and
$\vert H'\vert =\vert 5H_{\scriptscriptstyle{\sp 2}}
-2\sum_{\scriptscriptstyle {
1\leqslant i\leqslant 3}} x_i
-\sum_{\scriptscriptstyle {
4\leqslant j\leqslant 10}} y_j\vert$;
then $\pi ^{*}H$ is
very ample if and only if:
\smallskip 
\item{(0)} no two of the $y$'s are 
infinitely near each other;
at most one among the $y$'s 
is infinitely near to one of the $x$'s ;

\item{(1)} $h^0(H_{\scriptscriptstyle{\sp 2}}
-\sum_{i\in\Delta} x_i
-\sum_{j\in\Lambda} y_j)=0$, for 
$2\#\Delta +\#\Lambda\geqslant 5$ ;

\item{(2)} $h^0(2H_{\scriptscriptstyle{\sp 2}}
-\sum_{i\in\Delta} x_i
-\sum_{j\in\Lambda} y_j)=0$, for 
$2\#\Delta +\#\Lambda\geqslant 10$ ;

\item{(3)} $h^0(3H_{\scriptscriptstyle{\sp 2}}
-\sum_{\scriptscriptstyle {
1\leqslant i\leqslant 3}} x_i
-\sum_{\scriptscriptstyle {
4\leqslant j\leqslant 10}} y_j)=0$ .

\medskip

\noindent
(CF4) \ $deg\ W=11$ \ ($S$ is a hyperplane section of a 
Palatini scroll) 

\vskip -0.2 cm

\noindent
---------------------------------------------------------------------------------------- 

\noindent $W=\{ x_1,\ldots ,x_6,y_7,\ldots 
,y_{11}\}$, and
$\vert 6H_{\scriptscriptstyle {\sp 2}}-
\sum_{{\scriptscriptstyle{1\leqslant 
i\leqslant 6}}}2x_i-\sum_{{\scriptscriptstyle
{7\leqslant j\leqslant 11}}} y_j\vert ;$ 
then $\pi ^{*}H$ is
very ample if and only if:
\smallskip 
\item{(0)} at most one $y_j$ is
infinitely near to a point $x_i$ ;

\item{(1)} $h^0(H_{\scriptscriptstyle{\sp 2}}
-\sum_{i\in\Delta} x_i
-\sum_{j\in\Lambda} y_j)=0$, for 
$2\#\Delta +\#\Lambda\geqslant 6$ ;

\item{(2)} $h^0(2H_{\scriptscriptstyle{\sp 2}}
-\sum_{i\in\Delta} x_i
-\sum_{j\in\Lambda} y_j)=0$, for 
$2\#\Delta +\#\Lambda\geqslant 12$ ;

\item{(3)} $h^0(3H_{\scriptscriptstyle{\sp 2}}
-\sum_{\scriptscriptstyle {
1\leqslant i\leqslant 6}} x_i
-\sum_{\scriptscriptstyle{j\neq h}} y_j)=0$ 
for every $h\in\{ 7,\ldots ,11\} $ .


\bigskip\noindent
{\bf  2.-\ Birational projections of a threefold from a line on it.}
\smallskip


In the previous section we showed examples of
threefolds $X\subset\sp 5$ containing a line such that the projection 
from this line is a birational map $X \dashrightarrow \sp 3 .$
 Starting from this map it was very easy to find out a linear
system of surfaces in $\sp 3$ defining the birational inverse 
$\sp 3 \dashrightarrow  X$ of the projection. In this
section we want to formalize this procedure.

\smallskip

Let $X\subset\sp 5$ be a threefold satisfying our general 
assumptions; we will
denote by $d$ its degree. Moreover, we will assume that:

\smallskip
\noindent
{\sl there is a line $L$ on $X$  such that the projection 
 \ $\pi_L :X \dashrightarrow  \sp 3 $
from the center $L$  is birational.}

\tit {\bf Proposition 2.1} {\sl The surfaces in $\sp 3$ of the linear
system $\vert\Sigma\vert$ defining $\pi_L^{-1}$ have degree $d-1.$}

\proof This follows easily from the fact
 that any hyperplane $H\subset\sp 5$ intersects $L$ at a point.
\endproof

\smallskip

Let us denote by $M\simeq\sp 3$ the target of the projection $\pi_L.$
A key point here is the fact that, since $L$ and $M$ are disjoint, a
hyperplane $H$ containing $L$ cannot contain also $M,$ hence $H\cap M$
is a plane. From this remark we get at once the first part of

\tit {\bf Proposition 2.2} {\sl The 
set-theoretic image in $\pi_L$ of any hyperplane 
$H\subset\sp 5$
such that $L\subset H$ is a plane in $M.$  
The elements of $\vert\Sigma\vert$
corresponding to the hyperplanes in $\sp 5$ through $L$ break into a 
variable plane plus a fixed surface $\Phi ,$ of degree $d-2,$
which is the exceptional divisor for the map
$\pi_L^{-1}.$} 

\proof Let $u:\tilde X\longrightarrow X$ be 
the blowing-up of $X$ along $L,$ and
let $v:\tilde X\to\sp 3$ be the map which 
solves the singularities of 
$\pi_L :X \dashrightarrow \sp 3 .$ Then 
we have a commutative diagram

$$
\objectmargin {0.4pc}
\diagram
 \tilde X\dto_u \drto^v & \\
X \rdashed_{\pi_L}|>\tip   & \sp 3\cr
\enddiagram
$$

\smallskip

\noindent Let $\Bbb E$ denote the 
exceptional divisor of $u.$ Then, clearly,
a surface of $\vert\Sigma\vert$ corresponding 
to a hyperplane section  
$S\subset\sp 5$ through $L$ breaks 
into the plane $\pi_L(S)$ and
$\Phi :=v(\Bbb E ).$ So, by definition, 
the map $\pi_L^{-1}$ contracts $\Phi$ to $L.$
\endproof

\smallskip

\tit {\bf Remark 2.3} \ 
In certain cases we can use the particular 
structure of the surfaces in
$\vert\Sigma\vert$ corresponding to hyperplane 
sections through $L$ 
described above, to get informations on
the base locus of $\vert\Sigma\vert .$ Let 
$H\subset\sp 5$ be a hyperplane
such that $L\subset H;$ set $S=X\cap H$ and 
$V=\pi_L(S)\subset M.$ $V$ is a 
plane; the restriction of $\vert\Sigma\vert$ 
to $V$ is a linear system 
$\vert\Sigma '\vert$ of plane curves 
of degree $d-1.$ Of course,  
the dimension of $\vert\Sigma '\vert$ is 
$4,$ hence it defines a rational 
map whose image is a rational 
surface in $\sp 4,$ which has
the same degree $d$ as $X.$ The base 
locus of $\vert\Sigma '\vert$
is $V\cap B,$ where $B$ is the $1$-dimensional 
part of the base locus of
$\vert\Sigma\vert .$ In the case 
$\vert\Sigma '\vert$ is well understood and 
unique (see e.g. [A],[CF])
we get the desired informations about $B.$

\noindent
We will give an example of the use of 
this remark in the next section (Theorem 3.7).

\smallskip

\tit {\bf Proposition 2.4} \ 
{\sl The surface $\Phi$ is irreducible, 
rational and ruled. Every general plane 
section of $\Phi$ is a rational curve.} 

\proof In fact, since $\Bbb E$ is irreducible, 
$\Phi =v(\Bbb E )$ is also irreducible.

\noindent
If we think of the projection
$\pi_L$ as a rational map $\sp 5\backslash 
L\dashrightarrow M\simeq\sp 3,$ 
then $\Phi$ is the union of the images in 
$\pi_L$ of the (projectivized) 
tangent spaces to $X$ at the points of 
$L.$ In particular, from this it follows 
that $\Phi$ is ruled. In fact, let $A$ 
be a general point of $L,$ and denote by $T$
the tangent space to $X$ at $A.$ Then 
$\pi_L (T)$ is a line in $M,$ and
$\pi_L (T)\subset\Phi .$

The description of $\Phi$ by means of 
$\pi_L$ above shows
that the map $v$ restricted to $\Bbb E$ 
is still birational. 
Since $\Bbb E$ is a rational surface, 
$\Phi$ is also rational.

Finally, let $C$ be a general plane 
section of $\Phi.$ From
$\pi_L^{-1}(\Phi)=L$ it follows that 
$\pi_L^{-1}$ induces a 
birational map between $C$ and $L,$ 
hence $C$ is a rational curve.
\hfill\endproof

\tit {\bf Remark 2.5} As the examples 
treated in \S\thinspace 1
show, we may have on $\Phi$ more than 
one ruling of lines.
In any case, we agree to consider as 
{\sl the} ruling on $\Phi$
the one considered in the above proof.

\medskip

\tit {\bf Theorem 2.6} \ {\sl If there is a 
line $L\subset X$  
such that the projection  $\pi_L :X 
\dashrightarrow  \sp 3 $
from the center $L$  is birational, then $X$ is 
isomorphic to a blowing-down
of \ $\tilde{\Bbb P}^3,$ the blowing-up of 
$\sp 3$ along $B:=Bs(\vert
\Sigma \vert).$}

\proof
Let $u:\tilde X\longrightarrow X$ be the 
blowing-up of $X$ along $L$.
The homogeneous ideal of $L$ in $X$ is 
generated
in degree one; then the linear system  
$\vert 2H_X-L\vert$ on $X$
 defines the rational map
$ u^{-1} :X \dashrightarrow  \tilde X $ 
([Co],[BS]).
Now, the linear system $\vert H_X\vert$ 
on $\sp 3$ is 
$\vert (d-1)H_{\scriptscriptstyle {\sp 3}}
-B\vert .$ Finally, 
$\pi_L^{-1}$ contracts $\Phi$ to $L$ and 
$\Phi\in\vert (d-2)
H_{\scriptscriptstyle {\sp 3}}-B\vert .$ 
Therefore, on $\sp 3$
we have $\vert 2H_X-L\vert =\vert 
dH_{\scriptscriptstyle 
{\sp 3}}-B\vert .$
Let us recall that any smooth threefold 
of $\sp 5$
is linearly normal, hence the linear system
$\vert (d-1)H_{\scriptscriptstyle {\sp 3}}-
B\vert$ is complete.
Therefore the coherent sheaf 
$\hbox{\fascio I}_B(d-1)$ is spanned:
by [BS] it follows that
$\tilde X$ is isomorphic to the blow-up 
of $\sp 3$ along $B,$
namely, we have the commutative diagram

$$
\objectmargin {0.4pc}
\diagram
\tilde X\dto_u \rto^{\sim}& \tilde{\Bbb P}^3 
\dto\\
X  \rdashed_{\pi_L}|>\tip  & \sp 3\cr
\enddiagram
\leqno(1)
$$
\endproof

\smallskip

We remark that the linear system $\vert 
dH_{\scriptscriptstyle 
{\sp 3}}-B\vert$ defines a rational map
$ v:\sp 3 \dashrightarrow  \tilde X .$ This 
yields the factorization of $\pi _L^{-1}$

$$
\objectmargin {0.4pc}
\diagram
 & \tilde X\dto^u \\
\sp 3  \rdashed_{\pi_L^{-1}}|>\tip 
\urdashed^w |>\tip & X\cr
\enddiagram
$$

\smallskip

\noindent where $w$ is the birational 
inverse of the map 
$v:\tilde X\to\sp 3$ introduced in the 
proof of Proposition 2.2.

\medskip

To prove the next proposition we need the 
following

\smallskip

\noindent {\bf Lemma 2.7}\ {\sl Let $E\subset X$ 
be an irreducible exceptional 
divisor for the projection $\pi_L:X
\dashrightarrow\sp 3,$ and
set $B:=\pi_L(E).$ Assume that the base locus of
$\vert\Sigma\vert$ contains
exactly the $r$-th infinitesimal
neighbourhood of $B.$ If $dim(B)=1$ 
and if $b\in B$ is a general point, then
the fibre $\pi_L^{-1}(b)$ is a rational curve 
of degree $r+1,$ lying in the plane 
$\langle L\cup b\rangle .$ If $B$ is
a point, then $r=0$ and $\pi_L^{-1}(B)=
\langle L\cup B\rangle .$}

\proof Let $B$ be a curve.
It is clear that $\pi_L^{-1}(b)$ is a curve
lying in the plane $\langle L\cup b\rangle .$
The preimage $D\subset\tilde X$ of $b\in\sp 3$ 
in diagram $(1)$ is clearly a rational
curve in $\tilde X,$ which is not contained into the
exceptional divisor of $\tilde X.$ Therefore, the map
$u:D\to\pi_L^{-1}(b)$ is birational, hence $\pi_L^{-1}(b)$
is a rational curve. Finally, if $S$ denotes a
hyperplane section of $X,$ then $deg\thinspace \pi_L^{-1}(b)=
\pi_L^{-1}(b)\cdot S=r+1,$ since the base locus of
$\vert\Sigma\vert$ contains
exactly the $r$-th infinitesimal
neighbourhood of $B.$

Let $B$ be a point. Then, as a set, $\pi_L^{-1}(B)=
\langle L\cup B\rangle .$ Moreover, the general
curve section of $X$ intersects 
$\langle L\cup B\rangle$ transversally at a single
point. This shows that $r=0.$ 
\endproof

\tit {\bf Proposition 2.8} {\sl  \ 
Let $B$ be the base locus of $\vert\Sigma\vert .$ 
Then, we have $B\subset\Phi$ scheme-theoretically. 
If $R$ is a line for which the 
scheme-theoretic intersection $R\cap B$ is 
a zero dimensional scheme of length 
$d-1$ (\lq\lq\thinspace
$R$ is $(d-1)$-secant to $B$"),
then $R\subset\Phi$. If $B$ is
purely $1$-dimensional, then the converse is also
true (\lq\lq\thinspace the rational map 
$\pi^{-1}_L:B\dashrightarrow L$ has degree $d-1$")
and $H^0(\sp 3,\hbox{\fascio I}_B(d-2))$ is 
generated by $\Phi .$}

\proof The subset $\{ \Phi +H\ \big| \ H\in\check{\sp 3} \}
\subset\vert\Sigma\vert$ is still a linear system, whose base 
locus is exactly $\Phi .$ Therefore we have 
$\Phi\in H^0(\sp 3,\hbox{\fascio I}_B(d-2)).$

From $B\subset\Phi$ and from
the fact that the degree of $\Phi$ is $d-2,$ it
follows immediately that any line which is 
$(d-1)$-secant to $B$ lies on $\Phi .$ 

Conversely, assume that $B$ is purely 
$1$-dimensional and
let $R\subset\Phi$ be a line of the ruling. 
Therefore, there exists $\ell\in L$ such that 
$R=\pi_L(T_{\scriptscriptstyle{X,\ell}}).$
We set $X\cap T_{\scriptscriptstyle{X,\ell}}=L\cup D;$ then
$deg\ D=d-1.$ Assume that an irreducible component $D'$ of $D$
is not contained in the union of all exceptional divisors of
$\pi_L.$ Then $D'$ dominates $R$ through $\pi_L$ and the restriction
$D'\dashrightarrow R$ of $\pi_L$ is birational. This contradicts
$\pi_L^{-1}(R)=\ell .$ Hence, by Lemma 2.7, $D$ 
consists of finitely many
rational curves, each contained in a plane through $L,$ and the
sum of the degrees of these curves is $d-1.$ In view of Lemma
2.7 we conclude that the line $R$ is $(d-1)$-secant to 
the base locus $B.$

Finally, let $\Psi\in 
H^0(\sp 3,\hbox{\fascio I}_B(d-2)).$ 
Any $(d-1)$-secant line to $B$ is contained in $\Psi$
and since these lines fill $\Phi$ we conclude 
$\Phi\subseteq\Psi .$ But $\Phi$ and $\Psi$ have 
the same degree, hence $\Phi =\Psi .$
\endproof 

\medskip

The assumption that $\pi_L :X \dashrightarrow  \sp 3 $ is 
birational means that the general plane $U$ through $L$ 
cuts $X$ outside $L$ in
a single point. Let $U=H_1\cap H_2\cap H_3,$ where the $H_i$'s are hyperplanes
through $L.$ Now, $H_1\cap H_2\simeq\sp 3$ and we can think of $U$ as a plane 
inside this $\sp 3.$ Moreover, $X\cap H_1\cap H_2=L\cup\Delta.$ Therefore,  
$U$ cuts $\Delta$ outside $L$ in exactly one point.
Since the degree of $\Delta$ is $d-1,$ we conclude

\tit {\bf Proposition 2.9} {\sl The projection 
$\pi_L :X \dashrightarrow \sp 3 $  is birational if
and only if the general curve $\Delta$ defined above meets $L$ in $d
-2$ points.}

\smallskip

\noindent Moreover, concerning the curves $\Delta$ 
we can say a little bit more:

\smallskip

\noindent {\bf Proposition 2.10} {\sl The general 
curve $\Delta$ is rational.}

\proof The linear system on $X$ which defines the projection
$\pi_L :X \dashrightarrow  \sp 3 $ is \ 
$\vert H_X-L\vert .$
Now, for any $H_1,H_2\in\vert H_X-L\vert ,$ the curve $\Delta$ is the
\lq\lq\thinspace variable part" of $H_1\cap H_2.$ Therefore, $\pi_L(\Delta )$
is the line, intersection of the two planes  $H_1\cap M$ and $H_2\cap M.$
\endproof

\noindent This can be refined as follows.

\smallskip

\noindent {\bf Proposition 2.11} {\sl Let $\pi$ denote the sectional genus of
$X.$ Then the projection $\pi_L$ is birational if and only if for the 
general curve $\Delta$ it holds $\ga (\Delta )=\pi -d+3.$ In particular,
if $X$ has non special, rational hyperplane sections, then
$\pi_L$ is birational if and only if $\ga (\Delta )=0.$}

\proof We have the relation

$$
\pi =\ga (L\cup\Delta )=\ga (L)+\ga (\Delta )+deg \ (L\cap\Delta )-1
$$

\smallskip
\noindent and the first part of the statement follows from Proposition 2.6.
The last part follows from the fact that for such threefolds $\pi =d-3.$
\endproof
\smallskip
We conclude this section by reversing the point of view.

\tit {\bf Theorem 2.12} {\sl Let $B\subset\sp 3$ 
be a closed subscheme. We assume that 
for some positive integer $d:$ 
\item{(i)}
$\hbox{\fascio I}_B(d-1)$ 
is spanned  and 
$h^0(\hbox{\fascio I}_B(d-1))=6;$ 

\item{(ii)}$h^0(\hbox
{\fascio I}_B(d-2))=1;$

\item{(iii)}for some plane 
$H\subset\sp 3,$ \ $B\cap H$
is the base locus of a linear system of curves 
in $H,$ of
degree $d-1,$ which is very ample on the
blow-up of $H$ in $B\cap H.$

\noindent Then, denoted by
$X$ the image of the rational map $f:\sp 3
\dashrightarrow\sp 5$ defined by the linear 
system $\vert (d-1)H_{\scriptscriptstyle{\sp 3}}
-B\vert,$ it follows  
that $X$ is a smooth threefold of degree $d,$ 
the map $f$  is birational onto the image and
there is a line $L\subset X$ such that 
$f^{-1}$ is given by the projection from $L.$}

\proof By $(ii)$ we have $H^0(\hbox{\fascio I}
_B(d-2))=k\cdot\Phi$ and the linear system
$\vert (d-1)H_{\scriptscriptstyle{\sp 3}}
-B\vert$ contains the linear sub-system 
$\{ \Phi +H\ \big| \ H\in\check{\sp 3} \}$
of dimension $3$. Therefore, $f(\Phi )$ 
is contained in any hyperplane of a $3$-dimensional
linear family, hence  $L:=f(\Phi )$ is a line on $X.$

The restriction $\sp 3\backslash\Phi\to X\backslash L$
of $f$ is clearly a biregular map whose inverse
is the projection from $L.$ In particular, $deg \ X=d$
and $X$ is smooth outside $L.$ Finally,
by $(iii)$ there is a smooth hyperplane section
of $X$ containing $L,$ hence $X$ is smooth also along
$L.$
\endproof


\bigskip\noindent
{\bf  3.-\ The Castelnuovo threefold.}
\tit 


Let $X\subset\sp 5$ be a Castelnuovo threefold; 
it has degree $5$ and 
sectional genus $2.$  $X$ can be defined by the 
$2\times 2$ minors of a 
suitable $2\times 3$ matrix of forms; then it 
is arithmetically Cohen-Macaulay. The hyperplane 
sections of $X$ are Castelnuovo surfaces,
namely $S\subset\sp 4$ is the image of a rational 
map  
$ \phi :\sp 2 \dashrightarrow  \sp 4, $ 
birational onto $S$, defined
by the linear system $\vert 4H_{\scriptscriptstyle
{\sp 2}}-2x_1-\sum_{\scriptscriptstyle
{2\leqslant i\leqslant 8}} y_i\vert$ (see [A]). 
The images in $\phi$ of the lines through $x_1$ 
are conics on $S,$ which becomes in this way a
conic bundle over $\sp 1.$ It is easily seen 
that the bundle map is just the adjunction map. 
The adjunction map $\varphi :X\to \sp 1$ for $X$  
realizes the  Castelnuovo threefold as a quadric 
bundle over $\sp 1.$ (For more informations above 
the Castelnuovo threefold see [BOSS].)

\smallskip

In this section we will show the existence on any 
Castelnuovo 
threefold $X$ of a line $L$  such that 
$ \pi_L :X \dashrightarrow \sp 3$ is birational. 
Then we will apply the results of the previous
section to construct the linear system 
$\vert\Sigma\vert .$

\smallskip

\tit {\bf Proposition 3.1} {\sl Let 
$L\subset X$ be a line. The projection 
$\pi_L :X \dashrightarrow \sp 3$ is
birational if and only if $L$ is unisecant 
the quadrics on $X,$ i.e. $L$ is the image
of a section of $\varphi.$ In particular
$\varphi (L)=\sp 1.$}

\proof 
Assume that $\pi_L :X \dashrightarrow \sp 3$ 
is birational. Let $S\subset X$ be a hyperplane 
section containing $L.$ Since  the 
projection of $S$ from $L$ is still birational,
we have \  $\vert H_{\scriptscriptstyle {\sp 2}}\vert
\subset \vert 4H_{\scriptscriptstyle {\sp 2}}-2x_1-
\sum_i y_i\vert,$ hence there is an effective divisor 
$\Gamma\in\vert 3H_{\scriptscriptstyle {\sp 2}}
-2x_1-\sum_i y_i\vert ,$ whose image in $\phi 
:\sp 2 \dashrightarrow  S$ is the line $L.$

Note, in particular, that from this it follows that 
{\sl $L$ meets all the exceptional divisors on $S.$}

Now, since the cubic $\Gamma$ has a node in $x_1$ it 
intersect a general line through $x_1$ at exactly
one point away from $x_1.$ This is no more true for 
the lines $x_1y_i,$ $i>1,$ since $y_i\in\Gamma .$ 
But $\phi$ maps $x_1y_i$ onto a line on $S$ and the 
conic is given by this line together with the 
exceptional divisor over $y_i.$ As remarked above, 
$L$ intersects (in a unique point) also these conics.

\smallskip

Conversely, let $L\subset X$ be a line which is 
unisecant the quadrics on $X.$ Consider a pencil 
of hyperplanes whose center contains $L$ and the 
corresponding hyperplane sections they cut on $X.$ 
To this pencil it corresponds in the target
of the projection a pencil of planes 
with center a fixed line. Then the desired conclusion 
follows from next Lemma.
\endproof

\tit {\bf Lemma 3.2} {\sl Let $S$ be a Castelnuovo surface 
of $\sp 4$. If there exists a line $L$ contained in $S$
and unisecant its conics, then the linear system on $S$
of hyperplane sections containing $L$ 
defines a birational map to the plane.}

\proof  The condition \lq\lq\thinspace $S$ smooth" 
ensures that the base points $x_1, y_2\dots,y_8$ 
satisfy the condition (CF2), see 1.4.

$L$ being contained in a hyperplane section of $S$, 
it is the image under $\phi$ of a plane curve $D$ 
of degree at most $3$. It is easy to see 
that only three cases are possible:

\smallskip

(a) $D$ is a line passing
through three of the points $y_2\dots,y_8$, or

(b) $D$ is a conic passing through $x_1$ and five of 
the points $y_2\dots,y_8$, or

(c) $D\in \vert 3H_{\scriptscriptstyle {\sp 2}}
-2x_1-\sum_i y_i\vert.$

\smallskip

In all the three cases the residual system to $D$
is a homaloidal net in the plane; precisely, in 
case (a) it is the linear system of nodal cubics 
with node at $x_1$ and passing through $5$ fixed 
points, in case (b) it is the net of conics through
three fixed points, in case (c) it is the net of all
lines of the plane.
\endproof

\smallskip

The equations for $X$ are the $2\times 2$ minors of a matrix

$$
\hbox{\fascio M}=\left(\matrix{L_1 & L_2 & F_1\cr
L_3 & L_4 & F_2\cr}\right)
$$

\smallskip

\noindent where the $L_i$ are linear forms and the $F_j$ are quadratic
forms. 

\tit {\bf Lemma 3.3} {\sl Let $Q$ be the quadric hypersurface defined by 
$L_1L_4-L_2L_3=0;$ then $rk(Q)= 4.$}

\proof Clearly $rk(Q)\leqslant 4.$
Assume $rk(Q)\leqslant 3.$ Then $Q$ is a cone whose vertex
contains a plane $V;$ we can assume that $V$ is defined by the equations
$x_0=x_1=x_2=0.$ Then, after a change of coordinates, {\fascio M} has
the form

$$
\hbox{\fascio M}=\left(\matrix{x_0 & x_1 & F_1\cr
x_1 & x_2 & F_2\cr}\right).
$$

\smallskip

\noindent Therefore, by evaluating the jacobian matrix at a point $P\in V$ it
follows that the points $P\in V$ such that $F_1(P)=F_2(P)=0$
are singular for $X,$ a contradiction.
\endproof

\smallskip

\noindent Therefore $rk(Q)=4$ and the vertex of $Q$ is a line $L';$ 
assume that $x_0=x_1=x_2=x_3=0$ are equations for $L'.$ Then

$$
\hbox{\fascio M}=\left(\matrix{x_0 & x_1 & F_1\cr
x_2 & x_3 & F_2\cr}\right)
$$

\smallskip

\noindent Note, in particular, that $L'\subset X.$ The quadric $Q$ has 
equation $x_0x_3-x_1x_2=0.$ Then it is obtained by projecting from
$L'$ the smooth quadric defined by the same equation in some $\sp 3$
disjoint from $L'.$ Hence on $Q$ we have two $1$-dimensional families
of $\sp 3$'s \ $\big\{ M_{\mu}\big\}_{\mu\in\sp 1},$ \ 
$\big\{ N_{\nu}\big\}_{\nu\in\sp 1}$ \ such that

\smallskip

\item{-} each $M_{\mu}$ and $N_{\nu}$ contains $L';$ 

\item{-} $M_{\mu}\cap M_{\mu '}=L'$ \ and \ $N_{\nu}\cap N_{\nu '}=L';$

\item{-} $M_{\mu}\cap N_{\nu}$ \ is a plane.

\smallskip
\noindent{\bf Lemma 3.4} {\sl 
\item{(i)} The quadrics on $X$ are obtained as $G_{\mu}:=X\cap M_{\mu};$

\item{(ii)} $\big\{F_{\nu}:= X\cap N_{\nu}\big\}_{\nu}$ \ is a linear family of
cubic surfaces, each containing $L'.$

 }

\proof This follows easily from the equations for $X,$ since
\vskip -0.2 cm
$$
x_0-\mu x_2=x_1-\mu x_3=0\ \ \ \ \ \ \ \hbox{are equations for $M_{\mu}$}
$$
\vskip -0.2 cm
\noindent and
$$
x_0-\nu x_1=x_2-\nu x_3=0\ \ \ \ \ \ \ \hbox{are equations for $N_{\nu}$.}
$$
\vskip -0.2 cm
{}\endproof

\smallskip 
\noindent
{\bf Lemma 3.5} {\sl 
A general quadric $G_{\mu}$ \ and a general cubic
$F_{\nu}$ \ intersect in an irreducible conic.}

\proof We have $G_{\mu}\cap F_{\nu}=X\cap \big( M_{\mu}\cap N_{\nu}\big) .$
Therefore, by using equations it is easily seen that $G_{\mu}\cap F_{\nu}$
is the intersection of the plane $M_{\mu}\cap N_{\nu}$ with the irreducible
quadric of equation \ $\mu F_2-F_1=0.$
\endproof

\smallskip 
\noindent
{\bf Proposition 3.6} {\sl 
There exists on $X$ a line which is unisecant the quadrics.}

\proof We fix a general cubic surface 
$F_{\nu}$ and we denote $C_\mu :=
G_{\mu}\cap F_{\nu}$ the conics mentioned 
in the previous lemma; finally, we
recall that $L'\subset F_\nu.$
The surface $F_\nu$ cannot be singular
along $L'$ because otherwise
the conic $C_{\mu}$ would be reducible.
On the other hand, the base locus
of the linear system $\vert F_\nu \vert$
is the line $L'$, hence by Bertini
the general $F_\nu$ is smooth outside $L'.$
 From this it follows that on $F_{\nu}$ 
there exists a line
$L$ which is disjoint from $L'$ (see [BL]).

{\sl The line $L$ does not lie on any $M_{\mu}.$} 
Assume the contrary, and
let $L\subset M_{\mu}.$ Then $L$ is 
contained in the plane
$M_{\mu}\cap N_{\nu},$ which already 
contains $L',$ a contradiction to
$L\cap L'=\emptyset .$

Therefore, for every $\mu\in\sp 1$ the 
line $L$ intersects $M_{\mu}$
at a single point $U.$ Since $L\subset X,$ 
we have $U\in X\cap M_{\mu}=
Q_{\mu}$ and the proof is complete.
\endproof

\medskip

\noindent 
{\bf Theorem 3.7} {\sl Let $X$ be a Castelnuovo
threefold.

\item{(i)} There exists a line $L\subset 
X$ such that the projection $\pi _L$ of centre  
$L$ from $X$ to $\sp 3$ is birational. 

\item{(ii)} The birational map $\pi_L^{-1}$ 
 is defined by a linear system of quartic surfaces 
of $\sp 3$ whose base locus is the union of a 
hyperelliptic curve $B_2$ of degree $7$  
and arithmetic genus $3,$ having a \ $5$-secant line 
$B_1$, with the first infinitesimal neighbourhood 
of $B_1$. Moreover,
the exceptional divisors of  $\pi_L$  and 
  $\pi_L^{-1}$  are 
respectively 
$F_{\nu}\cup F'$, where $F'$ is a ruled 
surface of degree $9$, and 
 a rational 
ruled cubic surface $\Phi.$

}

\proof 
The first assertion follows directly
from 3.1 and 3.6.

To prove (ii) we may apply the results of $\S 2,$ 
in particular Remark 2.3. 
The linear system 
$\vert\Sigma\vert$ of 
surfaces of $\sp 3$ which defines $\pi_L^{-1}$ is
formed by quartics; the base locus $Bs(\vert\Sigma\vert)$
intersects a general 
plane in  $7$ simple points and one double point, 
or, to be more precise, the first infinitesimal neighbourhood
of a point. Therefore, $Bs(\vert\Sigma\vert)$ is the union
of the first infinitesimal neighbourhood of a 
line $B_1$ with a curve $B_2$ of degree $7$.
A characteristic curve $\Gamma$ of $\vert\Sigma\vert$
is a quintic of genus $2$.

To compute $p_a(B_2)$ and $deg \ (B_1\cap B_2)$ 
we will use the following degeneration argument.
Assume that $\vert\Sigma\vert$ contains two reducible 
surfaces of the form: $\Sigma _1=Q_1\cup Q_2$, $\Sigma_2
=M\cup D$ where $Q_1$ and $Q_2$ are quadrics, $M$ is a plane, 
$D$ is a cubic, all
containing $B_1$, and $Q_1\cap D$ contains a 
characteristic curve $\Gamma$. Then 
$$\Sigma _1\cap\Sigma _2=(Q_1\cap M)\cup (Q_1\cap D)\cup
(Q_2\cap M)\cup(Q_2\cap D)=$$
$$=(B_1\cup s)\cup (\Gamma\cup B_1)\cup (B_1\cup s')\cup 
(B_1\cup \Gamma ').$$
In this case $B_2$ splits as $B_2 =\Gamma '\cup s\cup
s'$. From this we compute $p_a(B_2)=3$ 
and $deg \ (B_1\cap B_2) =5$, 
moreover $deg \ (\Gamma\cap B_1)=3$ and 
$deg \ (\Gamma\cap B_2)=9$. 
In particular, $B_2$ is hyperelliptic, because it 
has degree $7$ and possesses the 
$5$-secant line $B_1$. Hence $B_2$ is contained 
in two cubic surfaces $G_1$, $G_2$ 
(see for instance [dA]), and is linked,
in the complete intersection of them, to a curve of degree
$2$ and genus $-2$: clearly this is a double structure on
$B_1$, therefore $G_1$ and $G_2$ are tangent along $B_1$.

The surface $\Phi$ is in this case a rational ruled 
cubic generated by the chords of
$B_2$ meeting $B_1$. The points of the double line
of the base locus come via $\pi _L$ from conics 
of planes containing $L$: they form a surface
 ruled by conics over $\sp 1$, whose degree is 
equal to $deg \ (\Gamma\cap B_1)=3$: 
it is precisely
the surface $F_{\nu}$ used in 3.6 to construct $L$.
The points of $B_2$ come from lines contained 
in $X$ and meeting $L$:
they form a ruled surface $F'$ of degree 
$9=deg \ (\Gamma\cap B_2)$, having $L$
as a double line.
\endproof

\smallskip
\noindent
{\bf Remark 3.8.} \ The lines on a Castelnuovo
threefold $X$ which are \lq\lq good centers of
projection" form a family of dimension $1.$

In fact, the lines we exhibit in the proof of Proposition
3.6 lie on a cubic surface $F_{\mu}\subset X.$ The 
surfaces $F_{\mu}$ are not ruled since they have at 
most isolated singularities. Hence, the dimension of
the family of lines we are interested in is at least $1.$

On the other hand, let us remark first that the only
quadrics on $X$ are the fibres of the adjunction map
$\varphi :X\to\sp 1$ (if $Q\subset X$ is a quadric
such that $\varphi (Q)=\sp 1,$ then $Q$ intersects
the general fibre of $\varphi$ in a line, and we get
a contradiction by a simple computation in $Pic(X)$).
A similar argument shows that if $V$ is a cubic
surface on $X,$ then $\langle V\rangle =\sp 3.$
Therefore, the hyperplanes of the pencil of center
$\langle G_{\mu}\rangle ,$ where $G_{\mu}$ is a
general quadric, intersect $X$ residually to $G_{\mu}$
exactly in the cubic surfaces $F_{\nu}.$

Now, if $L\subset X$ is a line which is unisecant
the quadrics on $X,$ for a fixed $G_{\mu}$ the linear
span $\langle L\cup G_{\mu}\rangle$ is a hyperplane
in $\sp 5$ and $L$ lies on the corresponding
cubic $F_{\nu}.$ 


\bigskip\noindent
{\bf  4.-\ Geometric properties of 
Bordiga and Palatini scrolls.}
\tit 


Let $X\subset\sp 5$ be a Bordiga scroll; it has degree $6$ and 
sectional genus $3.$  $X$ can be defined as the degeneracy locus of
a generic map of vector bundles $3\cO_{\scriptscriptstyle 
{\sp 5}}\to 4\cO_{\scriptscriptstyle {\sp 5}}(1);$ 
therefore it is arithmetically Cohen-Macaulay. 
A Bordiga scroll is a scroll over $\sp 2$,
more precisely $X=\Bbb P(\hbox{\fascio E})$
where {\fascio E} is a rank two bundle over 
$\sp 2$ with $c_1(\hbox{\fascio E})=4$. The scroll 
map $f :X\to \sp 2$ is the adjunction map.
The hyperplane sections of $X$ are Bordiga surfaces; 
in particular, $S\subset\sp 4$ is the image 
of a rational map (birational onto $S$) $ \phi 
:\sp 2 \dashrightarrow  \sp 4 $ defined by the 
linear system $\vert 4H_{\scriptscriptstyle {\sp 2}}
-\sum_{\scriptscriptstyle{1\leqslant i\leqslant 10}} 
y_i\vert$ ([A]). The inverse of $\phi$ is defined 
on the whole $S$ and is the adjunction map.

\smallskip

Let $X\subset\sp 5$ be a Palatini scroll; 
it has degree $7$ and sectional genus $4.$  
$X$ is arithmetically Buchsbaum; its coherent ideal
has a $\Omega$-resolution

\vskip -0.2 cm 

$$\objectmargin {0.4pc}
\diagram O \rto & 4\cO_{\scriptscriptstyle{\sp 5}}
 \rto^{\alpha} & \Omega^1_{\scriptscriptstyle{\sp 5}}
(2) \rto & \hbox{\fascio I}_X(4) \rto & 0 \cr
\enddiagram\leqno(1)
$$
\noindent where $\alpha$ is generic. 
A Palatini scroll is a scroll over a smooth cubic surface $V\subset\sp 3,$ 
and the scroll map $f :X\to V$ is still the adjunction map. 
If $S$ denotes
a general hyperplane section of $X$ we have the commutative diagram
$$
\objectmargin {0.4pc}
\diagram
X \rto^f & V \rto^g & \sp 2 \cr
S \umapsin \urto_{\psi} & & 
\enddiagram
\leqno(2)
$$
\noindent where $\psi :S\to V$ is the adjunction map 
for $S$ and $g:V\to\sp 2$ is the blow-up of $\sp 2$ 
at the points $x_1,\ldots ,x_6\in\sp 2.$
The map $g\circ\psi :S\to\sp 2$ is the blow-up of 
$\sp 2$ at  eleven points  $x_1,\ldots ,x_6,y_1,\ldots 
,y_5;$ more precisely, the linear system of curves in 
$\sp 2$ defining the birational map $(g\circ\psi)^{-1}$ 
is $\vert 6H_{\scriptscriptstyle {\sp 2}}-\sum_{
\scriptscriptstyle{1\leqslant i\leqslant 6}}2x_i
-\sum_{\scriptscriptstyle{1\leqslant j\leqslant 5}} 
y_j\vert$ (see [A]). 

\medskip

\noindent {\bf Proposition 4.1}\ {\sl If $S$ is a Bordiga surface then $S$ 
is a hyperplane section of a (smooth) Bordiga scroll of $\sp 5.$
If $S$ is a smooth, 
rational non special surface of degree $7$ in $\sp 4,$ 
then $S$ is a hyperplane section
of a unique arithmetically Buchsbaum threefold $X\subset\sp 5,$ whose
ideal has a $\Omega$-resolution like $(1).$ 
If $S$ is general, then
$X$ is smooth, i.e. it is a Palatini scroll.}

\proof (The following argument was inspired by [Ch].)

Both the Bordiga and Palatini scrolls are defined as degeneracy 
loci of a suitable map $\phi :\hbox{\fascio E}\to\hbox{\fascio F}$
of vector bundles over $\sp 5.$
Let $H\subset\sp 5$ be a general hyperplane.
From the exact sequence

\vskip -0.2 cm

$$
0 \to  {\check{\hbox{\fascio E}}}\otimes\hbox{\fascio F}(-1)
 \to  {\check{\hbox{\fascio E}}}\otimes\hbox{\fascio F} \to  
{\check{\hbox{\fascio E}}}\otimes\hbox{\fascio F}\vert_H \to  0
$$

\smallskip  

\noindent we get the long exact cohomology sequence

$$
H^0({\check{\hbox{\fascio E}}}\otimes\hbox{\fascio F}(-1))\hookrightarrow
H^0({\check{\hbox{\fascio E}}}\otimes\hbox{\fascio F}) \to 
H^0({\check{\hbox{\fascio E}}}\otimes\hbox{\fascio F}\vert_H) \to 
H^1({\check{\hbox{\fascio E}}}\otimes\hbox{\fascio F}(-1))
\leqno(3)
$$

\smallskip  

\noindent We recall that
$H^0({\check{\hbox{\fascio E}}}\otimes
\hbox{\fascio F})\simeq {\rm Hom}
(\hbox{\fascio E},
\hbox{\fascio F})$ and 
$H^0({\check{\hbox{\fascio E}}}
\otimes\hbox{\fascio F}\vert_H)\simeq 
{\rm Hom}(\hbox{\fascio E}\vert_H,\hbox{\fascio F}\vert_H).$

\smallskip

If $S$ is a Bordiga surface the last term in 
$(3)$ is zero.
Since $S$ is arithmetically Cohen-Macaulay, we conclude
that it is a hyperplane section of some smooth Bordiga scroll
$X$ ([HTV]). Note, however, that 
$h^0({\check{\hbox{\fascio E}}}\otimes\hbox{\fascio F}(-1)) 
=12,$ hence $X$ is not uniquely determined (this should be
compared with the case of the Palatini scroll below). 

\smallskip

If $S$ is a smooth, rational non special surface of 
degree $7$ in $\sp 4,$  by Bott'\thinspace s periodicity 
formulas, the first and the last term in $(3)$ are zero. 
Then any map $\overline\phi :
\hbox{\fascio E}\vert_H\to\hbox{\fascio F}\vert_H$ is the
restriction to $H$ of 
{\sl exactly one} map $\phi :\hbox{\fascio E}
\to\hbox{\fascio F}.$ This can also be seen in down to earth terms as 
follows. Any global section of $\Omega^1_
{\scriptscriptstyle {\sp n}}(2)$ can be identified with a
$(n+1)\times (n+1)$ skew-symmetric matrix with entries in 
the base field $k$ ([O]). Then $\phi :\hbox{\fascio E}
\to\hbox{\fascio F}$ corresponds to four $6\times 6$ skew-symmetric 
matrices. Let $x_5=0$ be the equation of $H$ in $\sp 5.$
Then a section of $\hbox{\fascio F}\vert_H$ corresponds to a $5\times 5$
skew-symmetric  matrix, plus a linear form $f=a_{05}x_0+a_{15}x_1+
\ldots +a_{45}x_4.$ Starting from $A$ and $f$ we can construct
the skew-symmetric $6\times 6$ matrix

\vskip -0.2 cm

$$\pmatrix{&&&a_{05}\cr
           &A&&\vdots\cr
           &&&a_{45}\cr
-a_{05}&\ldots &-a_{45}&0\cr},
$$

\noindent hence a section of $\Omega^1_
{\scriptscriptstyle {\sp 5}}(2).$ It 
is clear how to reverse this 
procedure and that this yields the 
identification ${\rm Hom}(\hbox{\fascio E},
\hbox{\fascio F})={\rm Hom}(\hbox{\fascio E}\vert_H,\hbox{\fascio F}\vert_H).$

A general map
$\overline\phi :\hbox{\fascio E}\vert_H\to\hbox{\fascio F}\vert_H$ is the
restriction to $H$ of a general map $\phi :\hbox{\fascio E}
\to\hbox{\fascio F},$ and we have the desired conclusion that $S$
is a hyperplane section of a (smooth) Palatini scroll.
\endproof
\medskip

We collect here three Lemmas, to be used in next Sections. 
The first one is due to F.L.Zak and S.L.L'vovsky, 
as well as the part of Lemma 4.3 concerning the
Palatini scroll ([ZL]).

\smallskip

\noindent {\bf Lemma 4.2}
{\sl Let $X\subset\sp r$ be an irreducible, non degenerate threefold
which is a scroll over the surface $V,$ and let $f :X\to V$ be
the scroll map. Let $Z\subset X$ be an irreducible surface such that
$dim\ \langle Z\rangle <r.$ Then either $f(Z)\subset V$ is a 
curve, or $f$ maps $Z$ birationally onto $V.$}

\proof If $f(Z)$ is not a curve, then $Z$ is mapped dominantly
onto $V.$ If this map is not birational, then the intersection of $Z$
with a fibre of $f$ is a zero dimensional scheme of length
greater than one. Hence the fibre of $f ,$ which is a line, is
contained in $\langle Z\rangle .$ This yields $X\subset\langle Z\rangle ,$
a contradiction since $X$ is non degenerate.
\endproof 

\smallskip

\noindent {\bf Lemma 4.3}
{\sl There are no planes on the scrolls of Bordiga and Palatini.
The Palatini scroll contains exactly $27$ quadrics; each of them
is smooth and has the form $f^{-1}(R),$ where $R$ is a line 
on $V,$ and $f$ denotes the scroll map.}

\proof  If $M\subset X$
is a plane, then obviously $f(M)$ cannot be a curve. By Lemma
4.2 we get in any case a birational map $f\vert_M:\sp 2\to\sp 2,$ 
whence $f\vert_M$ is an isomorphism. 
In the case of the Palatini scroll we reach immediately a contradiction 
because $f\vert_M$ factorizes through the blow-up $V\to\sp 2.$ 
In the case of the Bordiga scroll, the plane $M$ cuts on
a general hyperplane section $S\subset X$ one of the $10$ lines on $S$
(see Lemma 4.5 below). We have the commutative diagram

$$
\objectmargin {0.4pc}
\diagram
X \rto^f &  \sp 2 \cr
S \umapsin \urto_{\psi} &  \cr
\enddiagram
$$

\smallskip

\noindent where  $\psi :S\to\sp 2$ 
is the adjunction map for $S,$ which is the blow-up of $\sp 2$
at ten points. Therefore, any line on $S$
is contracted to a point, namely it is a fibre of $f .$
Therefore, $f$ contracts the plane $M$ to a point, 
a contradiction.

Let $X$ be a Palatini scroll, and
let $Q\subset X$ be an irreducible quadric. Since a quadric cannot
be mapped birationally onto a smooth cubic surface of $\sp 3,$ from
Lemma 4.2 it follows that $D=f(Q)$ is a rational curve on $V,$
and $Q=f^{-1}(D).$ 

Now, let {\fascio E} denote the $rk\ 2$ vector bundle 
on $V$ such that $X=\Bbb P (\hbox{{\fascio E}});$ we have 
$\cc (\hbox{{\fascio E}})=\cO_V(2)$ ([O]). Therefore, for any
curve $Y\subset V$: 

$$
deg\thinspace f^{-1}(Y)
=deg\thinspace\Bbb P({\hbox{\fascio E}\vert}_Y) =
\cc ({\hbox{\fascio E}})\cdot Y=2
\thinspace deg\thinspace Y.
$$

\smallskip

\noindent This shows that $D$ is a line and, conversely, 
if $D\subset V$ is a line, then
$f^{-1}(D)$ is a smooth quadric surface on $X.$ 
\endproof 

\medskip

Let $S$ be a hyperplane section of a Bordiga scroll, 
and let $u:S\to\sp 2$ be the blow-up whose inverse is 
defined by the linear system of plane curves $\vert 
4H_{\scriptscriptstyle{\sp 2}}-\sum_{\scriptscriptstyle
{1\leqslant i\leqslant 10}}x_i\vert .$
If $R\subset S$ is a line, assume $R=au^{*}H_
{\scriptscriptstyle {\sp 2}}-\sum_{\scriptscriptstyle
{1\leqslant i\leqslant 10}}b_iE_i$
in the Picard group of $S,$ where $E_i=u^{-1}(x_i),$ \ 
$a\geqslant 0$ and $b_i\geqslant 0$ for any $i.$ 
Since $R$ is properly contained in a hyperplane section
of $S,$ we can assume, moreover, that $a<4.$
Then $1=R\cdot H_S=4a-\sum_ib_i .$ If $a=1,$ 
three of the points $x_i$ should lie on a line, 
whereas if $a=2$ seven points $x_i$ should lie on a conic. 
Finally, $a=3$ implies $\sum_ib_i=11;$ but this is ruled out
by (CF3), since $S$ is smooth.

If $C\subset S$ is a conic, then $C$ is also contained 
in a hyperplane section $H_S$ of 
$S$ and it is different from $H_S.$ Then,
if $C=au^{*}H_{\scriptscriptstyle {\sp 2}}-\sum_
{i=1,\dots ,10}b_iE_i$ in $Pic(S),$ we can assume $0<a<4.$
From $2=C\cdot H_S=4a-\sum_ib_i$ it follows that for 
$a=2$ six of the points $x_i$ should lie on a conic. 
For $a=3$ all the points $x_i$ should lie on a cubic,
which is impossible by (CF3). 

A similar argument  can be applied to determine the conics
lying on a general hyperplane section of a Palatini 
scroll. We conclude:

\smallskip

\noindent {\bf Lemma 4.4}{\sl 

\item{(i)} Let $S$ be a smooth hyperplane section of a 
Bordiga scroll. If $h^0(H_{\scriptscriptstyle {\sp 2}}-\sum_
{i\in\Delta}x_i)=0$ for $\#\Delta\geqslant 3$ and
$h^0(2H_{\scriptscriptstyle {\sp 2}}-\sum_
{i\in\Delta}x_i)=0$ for $\#\Delta\geqslant 7,$ then
the lines on $S$ are only the exceptional divisors of $u.$ 
Moreover, if $h^0(2H_{\scriptscriptstyle {\sp 2}}-\sum_
{i\in\Delta}x_i)=0$ for $\#\Delta\geqslant 6,$
then the conics on $S$ correspond exactly to the $45$ lines 
$x_ix_j\subset\sp 2.$

\item{(ii)} Let \ $S$ \ be a smooth hyperplane section 
of a Palatini scroll \ $X.$ \ Assume 
that the inverse 
of the blow-up $u:S\to\sp 2$
is defined by the linear system
$\vert 6H_{\scriptscriptstyle {\sp 2}}-
\sum_{{\scriptscriptstyle
{1\leqslant i\leqslant 6}}}2x_i-\sum_{{\scriptscriptstyle
{1\leqslant j\leqslant 5}}} y_j\vert .$ 
If $h^0(H_{\scriptscriptstyle {\sp 2}}-\sum_
{i\in\Delta}x_i-\sum_{j\in\Lambda}y_j)=0$ 
for 
$2\#\Delta +\#\Lambda\geqslant 4$ (where 
$\Lambda\neq\emptyset $), and if
$h^0(2H_{\scriptscriptstyle {\sp 2}}-\sum_
{i\in\Delta}x_i-\sum_{j\in\Lambda}y_j)=0$ 
for \ $2\#\Delta +\#\Lambda\geqslant 10,$ \
then the conics on  $S$ Êare only: 
the $6$ exceptional 
 divisors 
$u^{-1}(x_i),$ the strict transforms of the 
$15$ lines $x_ix_h,$ and the strict 
 transforms 
of the $6$ conics through five of
the points $x_i.$ In particular, on 
 $S$ we 
have exactly $27$ conics.

}
\endproof

\smallskip

By Lemma 4.1  the hypotheses of the above lemma 
are fulfilled  for a general hyperplane section
of a general Palatini scroll. Then, from Lemma 4.3 
it follows at once 

\smallskip

\noindent {\bf Corollary 4.5}
{\sl The only conics on a general Palatini scroll $X$ are those 
contained in the $27$ quadrics on $X.$ Hence they form a (reducible)
$3$-dimensional family of plane curves which do 
not invade  $X.$
Moreover, the planes of these conics do not invade  $\sp 5.$}


\bigskip\noindent
{\bf  5.-\ Existence of a \lq\lq\thinspace good"
line on the 
Bordiga  scroll.}
\tit 


\smallskip
In this section $X$ will denote a smooth 
Bordiga scroll and $\phi$ the adjunction map for $X$.
As in \S 3 for the Castelnuovo threefold,  we will show 
the existence
on  $X$ of a line $L$ such that
$ \pi_L :X \dashrightarrow \sp 3
$ is birational. Then we will 
construct the linear system $\vert\Sigma\vert .$

\smallskip

\noindent {\bf Proposition 5.1}
{\sl Let $L\subset X$ be a line. Let $S$ be
a smooth hyperplane section of $X$ containing $L$ and 
$\phi _S:S\to \sp 2$ be the adjunction map.
The projection
$ \pi_L :X \dashrightarrow \sp 3
$ is birational if and only if 
$\phi _S(L)$ is a conic $\gamma .$
In this case, $\gamma$ passes through
$7$ of the points $p_1,\dots ,p_{10}$ which 
are the images of
the ten lines contracted by $\phi _S.$}

\proof
By Proposition 2.9, $\pi _L$ is birational 
if and only if $deg\ (\Delta \cap L)
=4$, where $\Delta$ is the residual of $L$ in a general
curve section of $X$.
Since clearly $L$ cannot be a line of the scroll, there 
are only
two possibilities, i.e.
either $\phi _S(L)$ is a line and  $\phi _S(\Delta )$
is a cubic, or both $\phi _S(L)$ and $\phi _S(\Delta )$
are conics. Only in the second case $deg\ (\Delta \cap L)
=4.$ Last assertion follows by remarking that the degree of 
$\phi _S^*(aH-b_1p_1-\dots -b_{10}p_{10})$ is 
$4a-b_1-\dots -b_{10}.$
\endproof

\smallskip

If $H$ is a hyperplane of $\sp 5$, the adjuction
map for $X\cap H$ will be denoted by $\phi _H.$ 
Each $H$ determines $10$  points of $\sp 2$ (not necessarily 
distinct)
are determined: the images of
the ten lines contracted by $\phi _H.$ 
Next lemma gives a \lq\lq\thinspace bound"
on the set of hyperplanes for which
these points are not distinct.

\smallskip

\noindent {\bf Lemma 5.2}
{\sl Let $X$ be a Bordiga scroll and $\check X$
denote its dual variety. Then a general point of 
$\check X$ represents a hyperplane which is 
tangent to $X$ at a unique point. The hyperplanes
which are tangent to $X$ along a line form a 
subset of codimension at least $2$ 
in $\check X$.}

\proof Let us recall that, if 
$dim \ \check X=4-h$, then the contact locus 
of a general tangent hyperplane is a linear 
space of dimension $h.$ In our case 
$dim \ \check X=4$ (see [E]), hence the first
assertion follows.
Let now $\widetilde X\subset X\times\check X$ denote the conormal variety
of $X$, and let $p,q$ its projections to 
$X$ and $\check X$.
Denote by $A$ the subset of $\check X$ 
representing hyperplanes $H$ such that 
$p(q^{-1}(H))$ is a line, 
hence $dim \ q^{-1}(H)=1$.
The general fibres of the restriction of $q$:
$q^{-1}(A)\to A$ have dimension $1$, therefore 
if $dim \ A=3$, then $dim \ q^{-1}(A)=4$ and
$q^{-1}(A)= \widetilde X$: a contradiction.
\endproof
\smallskip

We introduce the regular map
$f:\check\sp 5\to S^{(10)}(\sp 2)$, from the 
dual of $\sp 5$ to the tenth symmetric power of 
the plane, which takes $H$ to the images of
the ten lines contracted by $\phi _H.$
Let $V$ denote the image of $f$: it is projective
irreducible of dimension $5$. Our aim 
is to show that
$V$ meets the codimension two 
irreducible subvariety $W$
of $S^{(10)}(\sp 2)$
parametrizing $10$-uples of points, $7$ of them lying
on a conic.

Let us consider the natural 
map $p: (\sp 2)^{10}\to  S^{(10)}(\sp 2)$ and set 
$V'=p^{-1}(V)$, $W'=p^{-1}(W).$

\medskip

\noindent {\bf Lemma 5.3}
{\sl With the above notation, $V'\cap W'
\neq\emptyset .$ Hence $dim \ V'\cap W'\geqslant 3$.}
\proof
We shall use the following 
result of Debarre ([D]):

\smallskip

{\sl Let $\Bbb P$ be a product of projective spaces:
$\Bbb P=\sp {n_1}\times\dots\times\sp {n_r}$, 
and $Z,Y$ be closed irreducible subvarieties 
of $\Bbb P$. For any subset $I=\lbrace i_1,\dots, 
i_s\rbrace\neq \emptyset$ of $\lbrace 1, \dots ,r\rbrace$, 
denote by $p_I$ the projection 
$\Bbb P\to \sp I:=\sp {n_{i_1}}\times\dots\times
\sp {n_{i_s}}$. If $dim \ p_I(Z)+dim \ p_I(Y)
\geqslant dim \ \sp I$ for any $I$, then 
$Z\cap Y\neq\emptyset .$}

\smallskip

We shall apply this result in the situation
$\Bbb P=(\sp 2)^{10}$ to a pair of 
varieties $Z,Y$ where $Z$ (resp. $Y$) 
is an irreducible componenty of $V'$ (resp. 
$W'$).
The Debarre's condition becomes: 
$codim \ p_I(Y)\leqslant dim \ p_I(Z).$
It is easy to see that $codim \ p_I(Y) 
\leqslant 2$ for any $I$, so
we have simply to  exclude the possibility that
$dim \ p_I(Z)\leqslant 1$. Otherwise a general fiber of 
$p_I:Z\to p_I(Z)$,
corresponding to  points $q_1,\dots ,q_7$,
would have dimension at least
$4$, i.e. there would be a $4$-dimensional
family of hyperplane sections of the 
threefold $X$ containing the $7$ lines contracted
to $q_1,\dots ,q_7$ by the adjunction map: this 
possibility is 
excluded because these 
lines are pairwise disjoint.
\endproof

\medskip

\noindent 
{\bf Theorem 5.4} {\sl Let $X$ be a Bordiga scroll.

\item{(i)} There exists a line $L\subset X$ such 
that the projection $\pi _L$
of centre  $L$ from $X$ 
 to $\sp 3$
is birational. 

\item{(ii)} The birational map 
$\pi_L^{-1}$  is defined by a linear system
of quintic surfaces 
 of $\sp 3$ whose
 base locus is the union of the
first infinitesimal neighbourhood 
of a  cubic curve $B_1$ of arithmetic genus $0$
with a curve $B_2$ of degree $7$ and 
 $\ga (B_2)=0.$ 
 Moreover,  we  have  $deg \ (B_1\cap B_2)=12.$
 The exceptional 
 divisors are: for $\pi_L$  
a rational scroll $F$
of degree 
$8$, unbalanced of 
 type  $(1,7)$,  
and a rational surface of degree  $8$;  
for  $\pi_L^{-1}$
a rational ruled 
 quartic surface $\Phi .$

}

\proof
By Lemma 5.3 it follows that the above introduced
varieties $V$ and $W$ have non-empty intersection
of dimension $\geqslant 3$. Lemma 5.2 allows
to exclude the possibility that the
intersection is formed only by 
$10$-uples of non-distinct points. Hence
a general $10$-uple in $V\cap W$ 
represents $10$ distinct points, $7$ lying on 
a conic $\gamma$; it comes via $f$ from a 
hyperplane $H$ such that $\phi _H^{-1}(\gamma )$
is a line $L$ as required in (i) by 
Proposition 5.1.

By the results of \S 2, $\vert \Sigma \vert$
is a linear system of quintic surfaces, having 
as base locus the first infinitesimal neighbourhood 
of a cubic curve $B_1$ and a curve $B_2$ of 
degree $7$. The base locus is contained 
in a rational surface $\Phi$ of degree $4$, with 
rational plane sections, which is the
exceptional divisor for $\pi _L^{-1}.$
A  general surface like this
 is the projection
of a rational normal surface of $\sp 5$ from a 
 line and has  a  cubic of arithmetic genus $0$ as 
singular locus ([C]). 
Let us consider the surface $F$ defined as
the union of the lines of the scroll $X$ intersecting 
$L$, i.e. $F=\phi ^{-1}(\phi (L))$. It is   a smooth rational 
scroll and $deg \ F=c_1(\hbox{\fascio E})deg \ \phi (L)=8.$
Since $L$ is a unisecant line of $F$, the surface 
$F$ is  of type $(1,7)$.
The projection centered in $L$ contracts $F$
to a rational curve of degree $7$, which 
is precisely $B_2$.
$B_1$ comes from the surface $F'$ generated
by the conics contained in $X$ and having $L$ 
as a chord.
If $\Gamma$ is a characteristic curve of $\vert
\Sigma\vert$, of degree $6$, the degree
of $F'$ is equal to the number of intersections
of $\Gamma$ with $B_1$, i.e. $8$ ($\Gamma\cdot
\Sigma =30=deg \ X+2\thinspace deg \ (\Gamma \cap B_1)+
deg \ (\Gamma \cap B_2)$).

Finally we get $deg \ (B_1\cap B_2)=12$ by an easy 
computation in the Picard group of $\Phi '$,
regardless $\Phi '$ is balanced or not.
\endproof

\smallskip
\noindent
{\bf Remark 5.5.} There is only a finite number of
lines on a Bordiga scroll $X$ 
which are \lq\lq good centers of projection". 

In fact, a hyperplane section
$S$ of $X$ contains at most one of such lines, because
the adjunction map $\varphi :S\to\sp 2$ takes any such line  
to a conic, and any of these conics must contain
seven of the points images of the lines contracted by
$\varphi .$ Since the hyperplanes in
$\sp 5$ containing a fixed line form a $3$-dimensional
family, the intersection $V'\cap W'$ in Lemma 5.3 is
purely of dimension $3,$ and each of its components
corresponds to a \lq\lq good center of projection"
on $X.$


\bigskip\noindent
{\bf  6.-\  Non - existence of a \lq\lq\thinspace good"
line on the Palatini  scroll.}
\tit 


The purpose of this section is to prove the following

\medskip

\noindent {\bf Theorem 6.1}
{\sl Let $X\subset\sp 5$ be a Palatini scroll. 
Then there does not exist any line $L\subset X$
such that the projection  $\pi_L:X\dashrightarrow\sp 3$ 
is  birational.}

\proof
We shall  argue by contradiction. 

From the results of \S\ 2, if such a line exists, 
then the surfaces of the linear system 
$\vert\Sigma\vert$ have degree $6.$
In particular, from Remark 2.3 and from (CF4), it follows that the $1$-dimensional part $B$
of $Bs(\vert\Sigma\vert )$ decomposes into a curve $B_1$ of degree $5$ and
the first infinitesimal neighbourhood of a curve $B_2$ of degree $6.$
Let $S\subset X$ be a hyperplane section 
containing $L,$ and let
$H\subset\sp 5$ be a hyperplane containing $L,$ with $\langle S\rangle\neq H.$
Then 
$S\cap H=L\cup\Delta ,$
where $\Delta$ is a rational curve of degree $6,$ meeting $L$ at $5$
points by Proposition 2.9. With the notations of 
diagram $(2)$ (\S 4)
we set $C:=g\circ\psi (L)$ and
$D:=g\circ\psi (\Delta);$ the curves $C$ and $D$ are both irreducible and
rational and $C\cdot D=5.$ Therefore, either $C$ or $D$ is a line.
Since $\Delta$ is the variable part of $S\cap H,$ in the case $deg \ C=5$ and 
$deg \ D=1$ all the base points $x_1,\ldots ,x_6,y_1,\ldots ,y_5$ lie on $C.$
In the  case $deg \ C=1$ and $deg \ D=5,$ in order that the 
curves $D$ form
a net, the double points $x_1,\ldots ,x_6$ 
must lie on $D$ and the line
$C$ contains $y_1,\ldots ,y_5$ 
(actually, both cases are possible by (CF4)).

In the case $deg \ C=5$ the curve $C$ has a node at any point $x_i,$ then
the strict transform of $C$ in $f,$ which is nothing but $\varphi (L),$
is bisecant to any exceptional divisor $E_i\subset V$ of $f.$ Therefore, 
{\sl the
line $L$ is bisecant to any of the six quadric surfaces $\varphi^{-1}(E_i),$
which are pairwise disjoint.}

In the case $deg \ C=1$ the curve $\varphi (L)$ is bisecant to any line $R_i$ 
on $V,$ 
which is the strict transform of one of the six conics on $\sp 2$  through $5$
of the points $x_1,\ldots ,x_6.$ Hence {\sl the
line $L$ is bisecant to any of the six quadric surfaces $\varphi^{-1}(R_i),$
which are pairwise disjoint.}

By Lemma 2.7, any point of the double curve $B_2$ comes 
from a conic on $X$ with two points in common 
with $L.$  Therefore, the curve $B_2$ contains the 
six lines $L_i$ and, by degree, $B_2$ actually coincides 
with the union of the six lines $L_i$.

Now, the ruled surface $\Phi$ introduced in Proposition 2.2 has degree $5.$ 
Since the first infinitesimal neighbourhood of $B_2$ is contained in $\Phi,$ 
we have that 
$\Phi$ contains the trisecant lines to $B_2.$ If we take three lines
among the $L_i$'s, there is a quadric surface $Q\subset \sp 3$ containing them.
Therefore $Q\subset\Phi ,$ a contradiction since $\Phi$ is irreducible by
Proposition 2.4.
\endproof


\tit
{\bf  7.-\ The inverse problem: sufficient conditions
on the base locus.}
\tit


In the previous sections  
we showed that on $\sp 2\times\sp 1$ and on the
Del Pezzo, Castelnuovo and Bordiga 
threefolds it is always possible to find 
 a line $L\subset X$ such that
the projection $\pi_L :X\dashrightarrow\sp 3$ is
birational, and we described the
base locus $B$ of the linear system defining
$\pi_L^{-1}.$

In this section we will reverse the point of view
by studying the sufficiency of the conditions 
on $B$. The case of
$\sp 2\times\sp 1$ is straightforward.
It will turn out that in all the other cases
curves $B$ with components of degree, arithmetic 
genus and
reciprocal position as in Proposition
1.3.1, Theorem 3.7 and Theorem 5.4 respectively, 
and satisfying 
condition $(iii)$ of Theorem 2.12 are the 
base locus of  a linear 
system of surfaces (of degree $3,4,5$ 
respectively) which define  birational 
maps with smooth image; in other words 
for these curves all the hypotheses of 
Theorem 2.12 are fulfilled. 

\medskip

\noindent {\bf Theorem 7.1.}
{\sl Let $B\subset\sp 3$ be a locally
Cohen-Macaulay curve verifying one of the 
following conditions:

\item{(DP)}
$B$ has degree $5,$ 
 $\ga (B)=2$ and the general plane
section $Z$ of $B$ satisfies (CF1);

\item{(C)} 
$B$ is the union of 
a  curve $B_2$ of degree $7$ 
and arithmetic genus $3,$  with the first 
infinitesimal 
neighbourhood of a line $B_1$ such that $deg 
(B_1\cap B_2)=5;$ the general plane
section $Z$ of $B_1\cup B_2$ satisfies (CF2);

\item{(B)}
$B$ is the union of the
first infinitesimal neighbourhood
of a  cubic curve $B_1$ of arithmetic 
genus $0$
with a curve $B_2$ of degree $7$,
such that $\ga (B_2)=0$, 
and $deg (B_1\cap B_2)=12;$ the general plane
section $Z$ of $B_1\cup B_2$ satisfies (CF3').

\noindent
Then  $B$ fulfills the assumptions 
of Theorem 2.12, with $d=4,5,6$ respectively.}

\proof
The main lines of the proof are the same in all 
cases. More precisely, we investigate
the numerical character of $B$ (see [GP] for
the definition and first properties of the 
numerical character; see also [B]) and 
we prove that conditions (CF) determine it
uniquely. Using the properties of
this character we show that $B$ is arithmetically
Cohen-Macaulay and that $H^1({\cal O}_B(d-2))=0.$
Finally, we show that $\hbox{\fascio I}_B(d-1)$ is 
spanned by a liaison argument.
We shall give the details only for the 
Bordiga scroll.

\smallskip

To start, let us compute the arithmetic genus of $B.$
We denote the first infinitesimal neighbourhood of
$B_1$ by $2B_1.$ The arithmetic genus of $2B_1$ is
$8$; moreover, from $deg (B_1\cap B_2)=12$ it
follows $deg (2B_1\cap B_2)=24.$ Therefore, we have
$\ga (B)=31.$

Now, let $\Gamma$ be a general plane section of 
$B$ and $\sigma$ the minimum degree of a
 plane curve containing $\Gamma .$ First
of all, we observe that $\sigma =4.$ In fact:
$deg \ \Gamma =16,$ hence $\sigma\leqslant 5;$
it is easily seen that there is only one possibility
for the numerical character of $B$ if $\sigma =5,$
namely: $\chi =(6,5,5,5,5)$. But the genus of this
character is $g(\chi )=30$ which contradicts
$31=\ga (B)\leqslant g(\chi ).$ On the other hand,
$\sigma <4$ is impossible because, otherwise, $Z$ 
would be contained in
a plane cubic, hence contradicting (CF3').

For $\sigma =4$ the possibilities for the numerical 
character $\chi =(n_0,n_1,n_2,n_3)$ 
of $B$ are collected in the following table:

$$
\matrix{
& n_0 & n_1 & n_2 & n_3 &&\ \ \ \ \ \ \ \ && n_0 & n_1 & n_2 & n_3\cr
&&&&&&\ \ \ \ \ \ \ \ &&&&& \cr
a) & 10-r & 4+r  & 4  & 4, & 0\leqslant r\leqslant 2&\ \ \ \ \ \ \ \ &  e) &  6 & 6  & 6  & 4\cr
&&&&&&\ \ \ \ \ \ \ \ &&&&& \cr
b) &  7 & 7  & 4  & 4 & &\ \ \ \ \ \ \ \ &  f) &  7 & 5  & 5  & 5 \cr
&&&&&&\ \ \ \ \ \ \ \ &&&&& \cr
c) &  8 & 5  & 5  & 4 & &\ \ \ \ \ \ \ \ & g) &  6 & 6  & 5  & 5\cr 
&&&&&&\ \ \ \ \ \ \ \ &&&&& \cr
d) &  7 & 6  & 5  & 4 & & &&&& \cr}
$$

\smallskip

In the cases $a), b), c), e), f)$ the character 
is disconnected:  this implies that $B$ 
is the union of two subcurves $C_1$ and $C_2$
whose characters depend on the position 
of the first gap. In the cases $a), c), f)$ 
 the first gap is after $n_0$: 
hence one of the subcurves is plane, of degree respectively
$10-r, 8,7;$ this contradicts (CF3').

In the case $b)$  the gap is after $n_1,$ then one
of the subcurves is contained in a quadric; moreover 
we compute its degree to be $13:$ this 
also contradicts (CF3').

In the case $e)$ the gap is after $n_2;$ 
this implies that one 
of the two subcurves, say $C_1,$ is a line, while $C_2$
is contained in a cubic surface and has degree $15.$
If $C_1\subset B_1$
then $\Gamma$, and therefore also
$Z,$ is on a cubic, against (CF3'). If 
$C_1\subset B_2,$ then $C_2\supset B_1$ and
a general plane section 
of $C_2$ is contained in a plane cubic $D$, which 
necessarily has 
$3$ singular points and therefore splits; by 
considering the various possibilities 
for the splitting type of $D$,
we again get a contradiction to (CF3').

So the character $\chi$ is connected: 
we want to exclude
the possibility $d).$ In fact, if
$\chi =(7,6,5,4),$ then $\Gamma$ is contained in two
quartics $F, G.$ Since $B$ (and therefore also $\Gamma$)
is certainly not a complete intersection, by genus reasons,
we conclude that $F$ and $G$ have a common component,
of degree $h\leqslant 3$.
By [B], Proposition 1.5, $B$ 
is the union of two subcurves $C_1$ and $C_2$, 
and precisely:
if $h=3$, then one of them is a line, and we may argue
as in case $f)$; if $h=2$, then one of them has character
$(7,6)$ and has therefore degree $12$ and is contained
in a quadric; if $h=1$, then 
there is a plane component of degree $7$: in both
cases we reach a contradiction.

We have proved that $\chi =(6,6,5,5)$. Now: the 
speciality of $B$
is $\leqslant n_{\sigma -1}-2$, so $h^1({\cal O}_B(4))
=0$. Moreover $g(\chi )=31=p_a(B)$, which implies that $B$ 
is arithmetically Cohen-Macaulay. We get 
$h^0(\hbox{\fascio I}_B(4))=1$ and $h^0(\hbox{\fascio I}_B(5))=6$. 
We may perform  a liaison of type 
$(4,5)$, and find that the linked curve has degree $4$
and arithmetic genus $1$ and is therefore a complete 
intersection of two quadrics; by mapping cone we get
the minimal free resolution of $\hbox{\fascio I}_B$, 
showing that $\hbox{\fascio I}_B(5)$ is globally generated.
\endproof


\centerline{\bf References}


\medskip

\noindent\item {[A]} J. Alexander: 
\lq\lq\thinspace Surfaces 
rationnelles non-sp\'eciales dans $\sp 4$'',
 Math. Z. {\bf 200} (1988), 87-110

\vskip 0.1 cm 

\noindent\item {[A1]} J. Alexander: 
\lq\lq\thinspace Speciality one rational surfaces in
 $\sp 4$'', in Complex Projective Geometry 
(G.Ellingsrud, C.Peskine, G.Sacchiero, S.A.Stromme Eds.),
 Cambridge University Press (1992), 1-23

\vskip 0.1 cm 

\noindent\item {[dA]} J. D'Almeida: 
\lq\lq\thinspace Courbes de l'espace projectif: 
S\'eries lin\'eaires incompl\`etes et multis\'ecantes'',  
J. reine angew. Math. {\bf 370} (1986), 30-51

\vskip 0.1 cm

\noindent\item {[B]} V. Beorchia: 
\lq\lq\thinspace On the arithmetic genus of locally
Cohen-Macaulay space curves'',
Intern. J. of Math. {\bf 6} (1995), 491-502

\vskip 0.1 cm 

\noindent\item {[BOSS]} R. Braun - G. Ottaviani - M. Schneider -
F.-O. Schreyer: 
\lq\lq\thinspace Classification of conic bundles
in $\sp 5$'', Ann. Scuola Norm. Sup. Pisa {\bf 23} (1996),  69-97

\vskip 0.1 cm 

\noindent\item {[BL]}   M. Brundu - A. Logar: 
\lq\lq\thinspace Classification of cubic surfaces 
with computational methods'',
Quaderno Matematico n.375, Universit\` a di Trieste, (1996)

\vskip 0.1 cm 

\noindent\item {[BS]} M. C. Beltrametti - A. J. Sommese:
\lq\lq\thinspace Notes on embeddings of blowups'', 
preprint 1996, to appear in J. of Algebra

\vskip 0.1 cm 

\noindent\item {[CF]} F. Catanese - M. Franciosi:
\lq\lq\thinspace Divisors of small genus on surfaces and 
projective embeddings'', in
Israel Mathematical Conference Proceedings, Vol.9
(1996), 109-140
\vskip 0.1 cm

\noindent\item {[CH]} F. Catanese - K.Hulek: 
\lq\lq\thinspace Rational surfaces
in $\sp 4$ containing a plane curve'', preprint (1995)

\vskip 0.1 cm

\noindent\item {[Ch]} M. C. Chang: 
\lq\lq\thinspace On the hyperplane sections of certain
codimension $2$ subvarieties in $\sp n$'',  
Arch. Math. {\bf 58} (1992), 547-550

\vskip 0.1 cm

\noindent\item {[C]} F. Conforto: 
\lq\lq\thinspace {\sl Le superficie razionali''},
Zanichelli, Bologna, 1939

\vskip 0.1 cm

\noindent\item {[Co]} M.  Coppens: 
\lq\lq\thinspace Embeddings of 
blowing-ups'', in {\it Seminari di  Geometria 1991-93}, Universit\`a
di Bologna, 89-100  

\vskip 0.1 cm

\noindent\item {[D]} O.  Debarre: 
\lq\lq\thinspace Th\' eor\` emes de connexit\' e
pour les produits d'espaces projectifs et les
Grassmanniennes'', Amer. J. Math. {\bf 118} (1996), 1347-1367

\vskip 0.1 cm

\noindent\item {[E]} L.  Ein: 
\lq\lq\thinspace Varieties with small dual varieties, I'', 
Invent. math. {\bf 86} (1986), 63-74

\vskip 0.1 cm

\noindent\item {[EGA]}  A. Grothendieck - 
J. Dieudonn\' e: 
{\it \'El\' ements de G\' eom\' etrie Alg\' ebrique}
Chap. IV (Quatri\` eme Partie), 
Publ. Math. I.H.E.S. {\bf 32} 1967

\vskip 0.1 cm

\noindent\item {[GP]} L.  Gruson - Ch. Peskine: 
\lq\lq\thinspace Genre des courbes de l'espace projectif", 
in Algebraic Geometry, Troms\o  1977, LNM {\bf 687} (1978),
31-60

\vskip 0.1 cm

\noindent\item {[HTV]} J. Herzog - N.V. Trung - G. Valla: 
\lq\lq\thinspace On hyperplane sections of reduced
irreducible varieties of low codimension", 
J. Math. Kyoto Univ. {\bf 34} (1994), 47-72

\vskip 0.1 cm

\noindent\item {[I]} P. Ionescu: 
\lq\lq\thinspace Embedded projective varieties of small
invariants, II'', 
Rev. Roumaine math. pures appl. {\bf 31} (1986), 539-544

\vskip 0.1 cm

\noindent\item {[I1]} P. Ionescu: 
\lq\lq\thinspace Generalized adjunction and applications", 
Math.Proc.Camb. Phil.Soc. {\bf 99} (1986), 457-472

\vskip 0.1 cm

\noindent\item {[J]}   F. Jongmans:
\lq\lq\thinspace Les vari\'et\'es alg\'ebriques \` a trois 
dimensions dont les 
courbes - sections ont le genre trois'', Acad. Roy. Belgique,
Bull. Cl. Sci. (5), {\bf 30} (1943), 766-782, 823-835
 
\vskip 0.1 cm

\noindent\item {[M]}   U. Morin: 
\lq\lq\thinspace Sui tipi di sistemi lineari di superficie
algebriche a curva - caratteristica di genere due'', Ann. 
Mat. Pura Appl., Ser. IV {\bf 19} (1940), 257-288
 
\vskip 0.1 cm

\noindent\item {[M1]}   U. Morin: 
\lq\lq\thinspace Sulle variet\` a algebriche a 
curve - sezioni di genere tre'', Ann. 
Mat. Pura Appl., Ser. IV {\bf 21} (1942), 1-43

\vskip 0.1 cm

\noindent\item {[MP]} E.  Mezzetti - D. Portelli: 
\lq\lq\thinspace Linear systems representing threefolds
which are scrolls on a rational surface"
in preparation
 
\vskip 0.1 cm 

\noindent\item {[O]}   G. Ottaviani:
\lq\lq\thinspace On $3$-folds in $\sp 5$ which are scrolls'',
Ann. Scuola Norm. Sup. Pisa, Cl.
Sci. Ser. (4) {\bf 19} (1992),  451-471

\vskip 0.1 cm   

\noindent\item {[SR]}   J. G. Semple - L. Roth:\ 
{\it Introduction to Algebraic Geometry},
Clarendon Press, Oxford, 1949

\vskip 0.1 cm

\noindent\item {[Z]}   O. Zariski:
\lq\lq\thinspace Foundations of a general theory of birational
correspondences'',
Trans. Amer. Math. Soc. {\bf 53} (1943),  490-542

\vskip 0.1 cm   

\noindent\item {[ZL]} F. L. Zak - S. M. L'vovsky:
\lq\lq\thinspace 
Around Palatini variety'',
unpublished notes 
\bye